\documentclass[
11pt,twoside, a4paper,
]{amsart}

\title[Extended Picard complexes and linear algebraic groups]
{           {\protect\hfill \normalfont \tiny
06/12/06 
            \\ \vspace{10pt}}
Extended Picard complexes and linear algebraic groups
}
\author{Mikhail Borovoi}
\address{Borovoi: Raymond and Beverly Sackler School of Mathematical Sciences,
Tel Aviv University, 69978 Tel Aviv, Israel}
\email{borovoi@post.tau.ac.il}
\thanks{Borovoi was partially supported by the Hermann  Minkowski Center for Geometry}

\author{Joost van Hamel}
\address{van Hamel: K.U. Leuven, Departement Wiskunde,
  Celestijnenlaan 200B, B-3001 Leuven (Heverlee), Belgium}
\email{vanhamel@member.ams.org}


\usepackage[arrow,matrix]{xy}  
\newcommand{\diagram}[1]{#1}   

\usepackage{mathptmx}
\usepackage{amscd}
\usepackage{amssymb}
\usepackage{eucal}
\usepackage{amsxtra}

\usepackage{verbatim}

\DeclareSymbolFont{rsfs}{U}{rsfs}{m}{n}
\DeclareSymbolFontAlphabet{\mathcal}{rsfs}

\DeclareFontEncoding{OT2}{}{} 
\DeclareTextFontCommand{\textcyr}{\fontencoding{OT2}\fontfamily{wncyr}\fontseries{m}\fontshape{n}\selectfont}
\newcommand{\Sha}{\textcyr{Sh}}

\theoremstyle{plain}

\newtheorem{theorem}{Theorem}[section]
\newtheorem{proposition}[theorem]{Proposition}
\newtheorem{lemma}[theorem]{Lemma}
\newtheorem{corollary}[theorem]{Corollary}

\newtheorem{conditional-result}[theorem]{Conditional Result}

\newtheorem{introtheorem}{Theorem}
\newtheorem{introproposition}[introtheorem]{Proposition}

\newtheorem{introcorollary}[introtheorem]{Corollary}

\newtheorem*{theorem*}{Theorem}
\newtheorem*{proposition*}{Proposition}
\newtheorem*{lemma*}{Lemma}
\newtheorem*{corollary*}{Corollary}
\newtheorem*{question*}{Question}
\newtheorem*{conjecture*}{Conjecture}
\newtheorem*{claim*}{Claim}

\newtheorem*{introtheorem*}{Theorem}
\newtheorem*{introproposition*}{Proposition}
\newtheorem*{introlemma*}{Lemma}
\newtheorem*{introcorollary*}{Corollary}

\theoremstyle{definition}

\newtheorem{definition}[theorem]{Definition}

\newtheorem{construction}[theorem]{Construction}
\newtheorem*{definition*}{Definition}
\newtheorem*{example*}{Example}

\newtheorem{subsec}[theorem]{}

\theoremstyle{remark}
\newtheorem{remark}[theorem]{Remark}
\newtheorem{Remarks}[theorem]{Remarks}

{\end{theoremlist}\end{Remarks}}

\newcounter{thlistitem}
\renewcommand{\thethlistitem}{\roman{thlistitem}}
\newenvironment{theoremlist}{
\begin{list}{\makebox[1.5em]{\hfill\textup{(\thethlistitem)}}}%
{\usecounter{thlistitem}
\setlength{\leftmargin}{0em}
\setlength{\rightmargin}{0pt}
\setlength{\labelwidth}{-1em}
\setlength{\labelsep}{.5em}
}}{\end{list}}

\newcounter{assumlistitem}
\renewcommand{\theassumlistitem}{\roman{assumlistitem}}
  {\begin{list}
      {\hfill\textup{(\theassumlistitem)}}
      {\usecounter{assumlist}
       \setlength{\leftmargin}{0pt}
        \setlength{\rightmargin}{0pt}
        \setlength{\labelwidth}{-1.5em}}}
  {\end{list}}

\newtheorem*{remark*}{Remark}
\newtheorem*{Remarks*}{Remarks}
\newenvironment{remarks*}{\begin{Remarks*}\nopagebreak[4]
\rule{1em}{0ex}\par 
\begin{theoremlist}}%
{\end{theoremlist}\end{Remarks*}}

\newcommand{\mathb}[1]{\mathbf{#1}}

\newcommand{\capbar}{\overline}

\newcommand{\sH}{\mathcal{H}}
\newcommand{\sK}{\mathcal{K}}

\newcommand{\sA}{\mathcal{A}}
\newcommand{\sD}{\mathcal{D}}

\newcommand{\Acat}{\mathcal{A}}

\newcommand{\fcolon}{\colon} 
\newcommand{\iso}{\simeq}
\newcommand{\isoto}{\overset{\sim}{\to}}

\newcommand{\into}{\hookrightarrow}
\newcommand{\onto}{\twoheadrightarrow}

\newcommand{\labelto}[1]{\xrightarrow{\makebox[1.5em]{\scriptsize ${#1}$}}}

\newcommand{\id}{\operatorname{id}}
\DeclareMathOperator{\coker}{coker}
\newcommand{\cl}{\operatorname{cl}}

\newcommand{\Zz}{{\mathb{Z}}}

\newcommand{\Cc}{{\mathb{C}}}

\DeclareMathOperator{\Spec}{Spec}

\DeclareMathOperator{\Hom}{Hom}
\DeclareMathOperator{\Ext}{Ext}

\DeclareMathOperator*{\tensor}{\otimes}

\newcommand{\Gm}{\mathb{G}_\mathrm{m}}

\def\Gak{\overline{\mathb{G}}_\mathrm{a}}
\def\unipotent{^\mathrm{u}}
\def\uu{^\mathrm{u}}
\def\reductive{^\mathrm{red}}
\def\red{^\mathrm{red}}
\def\tor{^{\mathrm{tor}}}
\def\sc{^{\mathrm{sc}}}
\def\sss{^{\mathrm{ss}}}
\def\semisimple{^{\mathrm{ss}}}

\def\Tsc{{T\sc}}

\newcommand{\etale}{{\textup{\'et}}} 

\newcommand{\ab}{{\textup{ab}}}

\DeclareMathOperator{\Gal}{Gal}

\DeclareMathOperator{\Div}{Div}
\DeclareMathOperator{\sDiv}{{\sD\mathit{iv}}}
\DeclareMathOperator{\Pic}{Pic}

\newcommand{\cycles}{\mathcal{Z}}

\DeclareMathOperator{\divisor}{div}

\DeclareMathOperator{\Br}{Br}
\newcommand{\Bra}{\Br_\mathrm{a}}
\DeclareMathOperator{\Chi}{\mathb{X}}

\DeclareMathOperator{\U}{U}

\DeclareMathOperator{\UPic}{UPic}
\DeclareMathOperator{\KDiv}{KDiv}
\DeclareMathOperator{\Ob}{Ob}

\DeclareMathOperator{\ob}{ob}

\newcommand{\et}{\etale}

\newcommand{\kbar}{{\overline{k}}}

\newcommand{\Xbar}{{\capbar{X}}}
\newcommand{\Ybar}{{\capbar{Y}}}
\newcommand{\Zbar}{{\capbar{Z}}}

\newcommand{\Gbar}{{\capbar{G}}}
\newcommand{\Hbar}{{\capbar{H}}}

\newcommand{\Tbar}{{\capbar{T}}}

\def\Gmbar{{\capbar{\mathbf{G}}_{\textup{m}}}}
\def\Gabar{{\capbar{\mathbf{G}}_{\textup{a}}}}

\def\id{{\textrm{id}}}

\def\Gsc{{G\sc}}
\def\Tsc{{T\sc}}

\def\et{{\text{\'et}}}

\newcommand{\updot}{^{\scriptscriptstyle{\bullet}}}
\newcommand{\pmatr}[1]{\begin{pmatrix}#1\end{pmatrix}}

\def\gg{{\mathfrak{g}}}
\def\Mbul{{M\updot}}
\def\Hombul{\Hom\updot}
\def\Ibul{{I\updot}}

\begin{document}

\begin{abstract}
For a smooth geometrically integral variety $X$ over a field $k$ of characteristic 0,
we introduce and investigate the extended Picard complex $\textup{UPic}(X)$.
It is a certain complex of Galois modules of length 2,
whose zeroth cohomology is $\overline{k}[X]^\times/\overline{k}^\times$
and whose first cohomology is $\textup{Pic}(\overline{X})$,
where $\overline{k}$ is a fixed algebraic closure of $k$ and
$\overline{X}$ is obtained from $X$ by extension of scalars to $\overline{k}$.
When $X$ is a $k$-torsor of a connected linear $k$-group $G$,
we compute $\textup{UPic}(X)=\textup{UPic}(G)$ (in the derived category)
in terms of the algebraic fundamental group $\pi_1(G)$.
As an application we compute the elementary obstruction
for such $X$.
\end{abstract}
\maketitle

\section*{Introduction}

Throughout the paper, $k$ denotes a field of characteristic~$0$ and
$\kbar$ is a fixed algebraic closure  of $k$.
By a $k$-variety we mean a geometrically integral $k$-variety.
If $X$ is a $k$-variety, we write $\Xbar$ for $X\times_k \kbar$.

Let $G$ be a connected reductive  $k$-group.
Let
\begin{equation*}
 \rho \fcolon G\sc \onto G\sss \into G
\end{equation*}
be Deligne's homomorphism, where $G\sss$ is the derived subgroup of
$G$ (it is semisimple) and $G\sc$ is the universal covering of $G\sss$
(it is simply connected).  Let $T$ be a maximal torus of $G$ (defined
over $k$) and let $T\sc := \rho^{-1}(T)$ be the corresponding maximal
torus of $G\sc$.  The 2-term complex of tori
\begin{equation*}
  T\sc \labelto{\rho} T
\end{equation*}
(with $T\sc$ in degree $-1$)
plays an important role in the study of the arithmetic of reductive groups.
For example, the Galois hypercohomology $H^i(k,T\sc\to T)$ of this complex
is the abelian Galois cohomology of $G$ (cf. \cite{Borovoi:Memoir}).
The corresponding Galois module
\begin{equation*}
  \pi_1(\Gbar):=\Chi_*(\Tbar)/\rho_* \Chi_*(\Tbar\sc)
\end{equation*}
(where $\Chi_*$ denotes the cocharacter group of a torus)
is the algebraic fundamental group of $\Gbar$
(\emph{loc.~cit.}).
The related  group of multiplicative type over $\Cc$ with holomorphic $\Gal(\kbar/k)$-action
\begin{equation*}
 Z(\hat{G}):= \Hom(\pi_1(\Gbar), \Cc^\times) =
\ker[\Chi^*(\Tbar) \tensor \Cc^\times \to \Chi^*(\Tbar\sc) \tensor \Cc^\times]
\end{equation*}
(where $\Chi^*$ denotes the character group)
is the center  of a connected Langlands dual group $\hat{G}$ for  $G$,
considered by Kottwitz \cite{Kot84}.

Clearly, the above constructions
rely on the linear algebraic group structure of $G$.
However we show in this paper  that in fact they are related
to a very natural
geometric/cohomological construction that works for
an arbitrary smooth geometrically integral $k$-variety $X$.
Namely, we consider the cone $\UPic(\Xbar)$ of the morphism
\begin{equation*}
  \Gm(\kbar) \to \tau_{\leq 1} R \Gamma(\Xbar, \Gm)
\end{equation*}
in the derived category of discrete Galois modules.
More explicitly, this cone is represented by the $2$-term complex
\begin{equation*}
  {\kbar(\Xbar)}^\times/\kbar^\times \to \Div(\Xbar)
\end{equation*}
(with ${\kbar(\Xbar)}^\times/\kbar^\times$ in degree 0),
where $\kbar(\Xbar)$ denotes the field of rational functions on $\Xbar$,
and $\Div(\Xbar)$ is the divisor group of $\Xbar$.
It follows from the definitions that the cohomology groups  $\sH^i$
of the complex $\UPic(\Xbar)$ vanish for $i\neq 0,1$, and
\begin{align*}
 \sH^0(\UPic(\Xbar)) &= U(\Xbar) := \kbar[\Xbar]^\times/\kbar^\times\\
   \sH^1(\UPic(\Xbar)) &=  \Pic(\Xbar)
\end{align*}
where $\kbar[\Xbar]$ is the ring of regular functions on $\Xbar$.
We see that  $\UPic(\Xbar)$ can be regarded as a 2-extension
of the Picard group $\Pic(\Xbar)$ by $U(\Xbar)$.
We shall call $\UPic(\Xbar)$  \emph{the extended Picard complex of} $X$.
The importance of the extended Picard complex lies in the fact that
$\UPic(\Xbar)$ contains more information than $U(\Gbar)$ and $\Pic(\Gbar)$
separately.

Let $G$ be an arbitrary connected linear $k$-group, not necessarily reductive.
We write $G\uu$ for the unipotent radical of $G$, and set $G\red=G/G\uu$
(it is reductive).
We define $\pi_1(\Gbar):=\pi_1(\Gbar\red)$.
This means the following.
Let
\begin{equation*}
 \rho \fcolon G\sc \onto G\sss \into G\red
\end{equation*}
be Deligne's homomorphism, where  $G\sss$ is the derived subgroup of $G\red$
and $G\sc$ is the universal covering of $G\sss$.
Let $T$ be a maximal torus of $G\red$
and let $T\sc := \rho^{-1}(T)$ be
the corresponding maximal torus of $G\sc$.
Then $\pi_1(\Gbar)=\Chi_*(\Tbar)/\rho_* \Chi_*(\Tbar\sc)$.

Consider the
{derived dual complex to $\pi_1(\Gbar)$}, which
by definition is given by
$$
\pi_1(\Gbar)^D=(\Chi^*(\Tbar)\to\Chi^*(\Tbar\sc))
\text{\quad (with $\Chi^*(\Tbar)$ in degree 0). }
$$
By Rosenlicht's lemma \cite{Rosenlicht} we have $\U(\Gbar)=\Chi^*(\Gbar)$.
By a formula of Voskresenski\u\i\ \cite{Vos},  Fossum--Iversen \cite{Fossum-Iversen}
and Popov \cite{Popov},
we have $\Pic(\Gbar)\simeq \Chi^*(\ker[\rho\fcolon\Gbar\sc\to\Gbar\red])$.
From these results one can easily obtain  that
$$
\sH^i(\UPic(\Gbar))\simeq \sH^i(\pi_1(\Gbar)^D) \text{ for } i=0,1.
$$

The central result of this paper is
that $\UPic(\Gbar)$ and $\pi_1(\Gbar)^D$ themselves are
isomorphic in the derived category.
\begin{introtheorem}
[Theorem~\ref{thm:main}]
\label{thm:intro-main}
For a connected  group $G$ over a field $k$ of characteristic~$0$,
there is a canonical isomorphism, functorial in $G$,
$$
\UPic(\Gbar)\isoto\pi_1(\Gbar)^D
$$
in the derived category of discrete Galois modules.
\end{introtheorem}

Both Rosenlicht's lemma and the vanishing of $U(\Gbar)$ and $\Pic(\Gbar)$
for a  semisimple simply connected group $G$ are used in
the proof.

We also prove a version of Theorem \ref{thm:intro-main} for torsors.

\begin{introproposition}
[Lemma \ref{lem:upic-torsor}(iii)]
\label{prop:intro-torsors}
Let $G$ be a connected  group over a field $k$ of characteristic~$0$,
and let $X$ be a $k$-torsor under $G$.
There is a canonical isomorphism, functorial in $G$ and $X$,
$$
\UPic(\Xbar)\isoto\UPic(\Gbar)
$$
in the derived category of discrete Galois modules.
\end{introproposition}

\begin{introcorollary}\label{cor:intro-torsors}
Let $G$ be a connected  group over a field $k$ of characteristic~$0$,
and let $X$ be a $k$-torsor under $G$.
There is a canonical isomorphism, functorial in $G$ and $X$,
$$
\UPic(\Xbar)\isoto\pi_1(\Gbar)^D
$$
in the derived category of discrete Galois modules.
\end{introcorollary}

This central result gives a good conceptual explanation of many
existing results in the literature concerning the striking relationship
between the arithmetic of a linear algebraic group $G$
and the Galois modules  $\Chi^*(\Gbar)$ and  $\Pic(\Gbar)$.

\subsection*{Picard group and Brauer group}
\begin{introproposition}
[Corollary~\ref{cor:pic-brauer-upic}(\ref{corpart:pic-upic})]
\label{prop:intro-Pic}
Let $X$ be a smooth geometrically integral variety over $k$.
Then there is a canonical injection
$$
\Pic(X)\into H^1(k,\UPic(\Xbar))
$$
which is an isomorphism if $X(k) \neq \emptyset$ or if $\Br(k) =0$.
\end{introproposition}

\begin{introcorollary}\label{cor:intro-Pic-G}
For a connected linear algebraic group $G$ over $k$
we have a canonical isomorphism
$$
\Pic(G)\isoto H^1(k,\pi_1(\Gbar)^D).
$$
\end{introcorollary}
\begin{proof}
The corollary follows immediately from
Proposition \ref{prop:intro-Pic} and Theorem \ref{thm:intro-main}.
\end{proof}

Let $X$ be a smooth geometrically integral variety over $k$.
Let $\Br(X) = H_{\etale}^2(X, \Gm)$ be the Brauer group of $X$, and let
$\Br_1(X)$ be the kernel of the map $\Br(X) \to \Br(\Xbar)$.
We write $\Bra(X)$ for the cokernel of the
canonical homomorphism $\Br(k)\to\Br_1(X)$.

\begin{introproposition}
[Corollary~\ref{cor:pic-brauer-upic}(\ref{corpart:brauer-upic})]
\label{prop:intro-Br}
 Let $X$ be a smooth geometrically integral variety over $k$.
There is a canonical injection
\begin{equation*}
  \Bra(X) \into H^2(k,\UPic(\Xbar) )
\end{equation*}
which is an isomorphism if $X(k) \neq \emptyset$ or $H^3(k,\kbar^\times) = 0$.
\end{introproposition}

\begin{introcorollary}\label{cor:intro-Br-G}
For a connected linear algebraic group over $k$ we have
a canonical  isomorphism
\begin{equation*}
  \Bra(G) \isoto H^2(k, \pi_1(\Gbar)^D)
\end{equation*}
\end{introcorollary}
\begin{proof}
  This follows immediately from Proposition \ref{prop:intro-Br} and
  Theorem \ref{thm:intro-main}.
\end{proof}
Note that Corollaries \ref{cor:intro-Pic-G} and \ref{cor:intro-Br-G}
are versions of results of Kottwitz~\cite[2.4]{Kot84}.  Kottwitz
proved that for a connected reductive $k$-group $G$ we have
\def\gg{{\mathfrak{g}}}
$$
\Pic(G)=\pi_0(Z(\hat{G})^\gg),\quad
 \Bra(G) =H^1(k,Z(\hat{G})),
$$
where $\gg=\Gal(\kbar/k)$.

\subsection*{UPic and smooth compactifications}

\begin{introproposition}[Proposition \ref{prop:zariski-open-upic}]\label{prop:intro-compactification}
Let $Y$ be a smooth compactification of a smooth geometrically integral $k$-variety $X$.
Then we have a distinguished triangle
$$
\Pic(\Ybar)[-1]\labelto{j^*}\UPic(\Xbar)\to \cycles^1_{Y-X}\to \Pic(\Ybar)
$$
where the morphism $j^*$ is induced by the inclusion map $j\fcolon X\to Y$,
and $\cycles^1_{Y - X}$ is the permutation module of
divisors in the complement of\ \ $\Xbar$ in $\Ybar$.
\end{introproposition}

We see that $\Pic(\Ybar)$ is very close to $\UPic(\Xbar)$:
up to translation, the difference between them is a permutation module.

If $C$  is a complex of $\Gal(\kbar/k)$-modules, we write
$$
\Sha^i_\omega(k, C)=\ker\left[H^i(k,C)\to \prod_\gamma H^i(\gamma,C)\right]
$$
where $H^i(k,C)$ is the corresponding Galois hypercohomology group, and
$\gamma$ runs over all closed procyclic subgroups of $\Gal(\kbar/k)$.

\begin{introproposition}[Corollary \ref{cor:h1-pic-open-dense}]\label{prop:intro-Sha}
Let $Y$ be a smooth compactification of a smooth $k$-variety $X$. Then
there is a canonical isomorphism
$$
\Sha^1_\omega(k,\Pic(\Ybar))\isoto \Sha^2_\omega(k,\UPic(\Xbar)).
$$
\end{introproposition}
Proposition \ref{prop:intro-Sha} follows easily
from Proposition \ref{prop:intro-compactification}.

\begin{introcorollary}
Let $Y$ be a smooth compactification of a $k$-torsor $X$ under a
connected linear  $k$-group $G$.
There is a canonical isomorphism
\begin{equation*}
\Sha^1_\omega(k, \Pic(\Ybar)) \simeq \Sha_\omega^2(k, \pi_1(\Gbar)^D)
\end{equation*}
\end{introcorollary}
\begin{proof}
The corollary follows immediately from Proposition \ref{prop:intro-Sha}
and Corollary \ref{cor:intro-torsors}.
\end{proof}
Note that we have $H^1(k,\Pic(\Ybar))=\Sha^1_\omega(k, \Pic(\Ybar))$
(see~\cite[Prop.~3.2]{CTK98}, \cite[Cor. 3.4]{Borovoi-Kunyavskii:arithmetical}).
Thus we have a new proof of the fact that
$$
H^1(k,\Pic(\Ybar))\simeq \Sha_\omega^2(k, \pi_1(\Gbar)^D),
$$
cf.~\cite[Thm. 2.4]{Borovoi-Kunyavskii:formulas}.

\subsection*{Elementary obstruction}

Let $X$ be a smooth geometrically integral $k$-variety.
We have an extension of complexes of Galois modules
$$
0\to\kbar^\times\to\left({\kbar(\Xbar)}^\times\to \Div(\Xbar)\right)\to
       \left({\kbar(\Xbar)}^\times/\kbar^\times\to \Div(\Xbar)\right)\to 0.
$$
It defines an element $e(X)\in\Ext^1(\UPic(\Xbar),\kbar^\times)$.  If
$X$ has a $k$-point, then this extension splits (in the derived
category), hence $e(X)=0$.
We shall call
$e(X)$ the \emph{elementary obstruction}
to the existence of a $k$-point in $X$, since it is a variant of the
original elementary obstruction of
Colliot-Th\'el\`ene and Sansuc \cite[D\'ef.\ 2.2.1]{CT-Sansuc}
which lives in $\Ext^1({\kbar(\Xbar)}^\times/\kbar^\times, \kbar^\times)$.

Now let $G$ be a connected linear $k$-group
and let $X$ be a $k$-torsor under $G$.
By Corollary \ref{cor:intro-torsors} we have  $\UPic(\Xbar)=\pi_1(\Gbar)^D$.
Using  Lemma \ref{lem:Ext-H^i} below, we obtain
$$
\Ext^1(\UPic(\Xbar),\kbar^\times)=H^1(k,\Hom(\pi_1(\Gbar)^D,\kbar^\times))
=H^1(k,T\sc\to T)
$$
(where $T\sc$ is in degree $-1$).
Recall that
the first abelian Galois cohomology group of $G$ is by definition
the abelian group $H^1_{\ab}(k,G):=H^1(k,T\sc\to T)$,
so the above identification gives us
$e(X)\in H^1_{\ab}(k,G)$.
Here we compare the elementary obstruction $e(X)\in H^1_{\ab}(k,G)$
with the image of the cohomology class $[X]\in H^1(k,G)$ of the torsor $X$
under the abelianization map $\ab^1\fcolon H^1(k,G)\to H^1_{\ab}(k,G)$
 constructed in \cite{Borovoi:Memoir}.

\begin{introtheorem}[Theorem \ref{th:elementary-abcoh}]
\label{thm:intro-torsor-obstruction}
Let $X$ be a $k$-torsor under a connected linear $k$-group $G$.
With notation as above, we have  $e(X)=\ab^1([X])$.
\end{introtheorem}

The theorem allows us to
translate existing results on abelian Galois cohomology of connected $k$-groups
to results on the elementary obstruction for torsors.
We simultaneously obtain results on smooth compactifications of
torsors,
since Proposition~\ref{prop:intro-compactification} implies
that the elementary obstruction $e(Y)$ for a smooth compactification $Y$
of a smooth variety $X$
vanishes if and only if the elementary obstruction $e(X)$ for $X$
vanishes.

\begin{introproposition}[Proposition \ref{prop:torsor-p-adic}]
For (a smooth compactification of) a torsor
under a connected linear algebraic
group $G$ over a $p$-adic field $k$, the
elementary obstruction is the only obstruction to the existence of
$k$-rational points.
\end{introproposition}

\begin{introproposition}[Proposition \ref{prop:torsor-number}]
\label{prop:intro-torsor-number}
For (a smooth compactification of) a torsor under a connected linear algebraic
group $G$ over a number field $k$, the
elementary obstruction is the only obstruction to the Hasse
principle.
\end{introproposition}

\begin{introcorollary}[Sansuc \cite{Sansuc:brauer-gal}, Cor. 8.7]
For a smooth compactification $Y$ of a torsor $X$
under a connected linear algebraic
group $G$ over a number field $k$, the
Brauer--Manin obstruction is the only obstruction to the Hasse
principle.
\end{introcorollary}
\begin{proof}
 Assume that $Y$ has points over all the completions of $k$.
  By \cite[Prop.~6.1.4]{Skorobogatov:torsors} the vanishing of the
  Brauer--Manin obstruction implies that the elementary obstruction
  vanishes, and we see from Proposition \ref{prop:intro-torsor-number} that $Y$
has a $k$-point.
\end{proof}

The results of this paper were announced in \cite{BvH06}.
\bigskip

{\bf Acknowledgements.}
The authors are very grateful to K.F. Lai for the invitation of M. Borovoi
to the University of Sydney, where their collaboration started, and to
J. Bernstein for most useful advice and for proving Lemma \ref{lem:Ext-H^i}.
We are grateful to V. Hinich and B. Kunyavski\u\i\ for useful discussions.
The first-named author worked on this paper while visiting the Max-Plank-Institut f\"ur
Matematik (Bonn) and Ohio State University; the hospitality and support of these institutions
are gratefully acknowledged.


\section{Preliminaries}\label{sec:preliminaries}

Throughout this paper, $k$ will be a field of characteristic zero.
Let $\kbar$ denote a fixed algebraic closure of $k$.  For a variety
$X$ over $k$ we denote by $D^{b}(X_\et)$ the derived category of
complexes of sheaves on the (small) \'etale site over $X$ with bounded
cohomology.  We write
\begin{equation*}
R \Gamma_{X/k} := R \varphi_* \fcolon D^{+}(X_\et) \to D^+(k_\et).
\end{equation*}
where $\varphi \fcolon X \to \Spec k$ denotes the structure morphism.
 We shall not distinguish between the category of \'etale
sheaves on $\Spec k$ and the category of discrete Galois modules.
We shall always assume our varieties to be geometrically integral.

Let $\Gm$ be the multiplicative group.  We shall
denote an \'etale sheaf represented by a group scheme by the same
symbol as the group scheme itself.
For a  variety $X$ over $k$ write $\Xbar=X\times_k \kbar$.
We define the following Galois modules:
\begin{align*}
  \U(\Xbar) & := (\Gamma_{X/k}\Gm)/\Gm =\kbar[\Xbar]^\times / \kbar^\times \\
  \Pic(\Xbar) & := R^1 \Gamma_{X/k} \Gm=H^1(\Xbar,\Gm).
  \end{align*}
These Galois modules are contravariantly functorial in $X$.

In this paper we shall be mostly interested in a \emph{complex} of
Galois modules that combines $\U(\Xbar)$ and $\Pic(\Xbar)$.  For this we want
to take the object $\tau_{\leq 1} R \Gamma_{X/k} \Gm$ in
$D^{b}(k_\et)$ modulo $\Gm$ (i.e. modulo $\kbar^\times$), where $ R
\Gamma_{X/k}$ is the derived functor, and $\tau_{\leq 1}$ is the
truncation functor.  To make this precise, we shall introduce some
terminology and notation.  For definitions of derived categories,
triangulated categories, derived functors, truncation functors etc. we
refer to original works \cite{Ver77}, \cite{Ver96}, \cite{BBD}, and
textbooks \cite{Iversen}, \cite{GM}, \cite{Weibel} (see also \cite{GM94}).

\begin{subsec}\label{subsec:cones-and-fibres}
{\it Cones and  fibres.}
Let $f\fcolon P \to Q$ be a morphism of complexes of objects of an
abelian category $\Acat$.
We denote by
$$\langle P \to Q ]$$
the \emph{cone} of $f$, i.e., the complex with the object in
degree $i$ equal to
$$P^{i+1} \oplus Q^{i},$$
and differential given by the matrix
\begin{equation*}
\begin{pmatrix}
-d_P & 0 \\
- f  & d_Q
\end{pmatrix},
\end{equation*}
which denotes the homomorphism
$(p, q) \mapsto (-d_P(p), - f(p)+ d_Q(q)  )$.
We adopt the convention that the diagrams of the form
\begin{equation}\label{eq:triangle-cone}
P \labelto{f} Q \labelto{\pmatr{0 \\ \id}} \langle P \to
Q]\labelto{\pmatr{\id & 0}} P[1]
\end{equation}
are distinguished triangles.

Similarly, we denote by
$$[ P \to Q \rangle$$
the complex with the object in degree $i$ equal to
 $$P^{i} \oplus Q^{i-1},$$
and differential given by the matrix
\begin{equation}\label{eq:matrix-fibre}
\begin{pmatrix}
d_P & 0 \\
 f  & -d_Q
\end{pmatrix},
\end{equation}
which denotes the homomorphism
$(p, q) \mapsto (d_P(p),  f(p)-d_Q(q)  )$.
We call $[ P \to Q \rangle$ the \emph{fibre} (or co-cone) of $f$.
Then
\begin{equation*}
\langle P \to Q ] = [ P \to Q \rangle [1],
\end{equation*}
and we have a distinguished triangle
\begin{equation}\label{eq:triangle-fibre}
 [ P \to Q \rangle \labelto{\pmatr{-\id & 0}} P \labelto{f} Q
\labelto{\pmatr{0 \\ \id} } [ P \to Q \rangle[1].
\end{equation}
We have
$[ P \to 0 \rangle = P$, and $\langle 0 \to Q ] = Q$.
\end{subsec}

\begin{remark}
Note that our sign convention for the differentials in the cone
corresponds to
\cite[I.4]{Iversen}, but is different from other sources,
such as~\cite[III.3.2]{GM}.
For example,
in the latter the cone has differential
\[
\begin{pmatrix}
-d_P & 0 \\
 f  & d_Q
\end{pmatrix}.
\]

The choice of signs also has an influence on the class of
distinguished triangles.
Indeed,
consider the following diagram
\begin{equation}\label{eq:GM-triangle-cone}
P \labelto{f} Q \labelto{\pmatr{0 \\ \id}}
C_\textup{GM}(f)\labelto{\pmatr{\id & 0}} P[1]
\end{equation}
where we write $C_\textup{GM}(f)$ for the cone as defined
in~\cite{GM}.
Then this diagram is a distinguished
triangle in $D(\Acat)$ in the convention of~\cite{GM} (cf. Def.~III.3.4 and
Lemma~III.3.3 in \emph{loc.~cit.}).
However, in our convention we would need to change the last
homomorphism of diagram~\eqref{eq:GM-triangle-cone}
to $\pmatr{-\id & 0}$ in order to have a distinguished
triangle.
\end{remark}

\begin{subsec}
Let $f\fcolon P\to Q$ be a morphism in the derived category $D^b(\sA)$.
We define a cone $\langle P\to Q]$ as the third vertex
of a distinguished triangle \eqref{eq:triangle-cone}.
Similarly, we define a fibre $[P\to Q\rangle$ as the third vertex
of a distinguished triangle \eqref{eq:triangle-fibre}.
It is well known that in general in a derived category
(or in a triangulated category) a cone
and a fibre are defined only up to a non-canonical isomorphism.
However we shall prove, that all the cones and fibres that we shall consider, will be
defined up to a \emph{canonical} isomorphism
(we shall use \cite[Prop. 1.1.9]{BBD}).
\end{subsec}

\begin{subsec}
  \textit{$\Ext$ and Galois cohomology.}  In order to compute the
  elementary obstruction to the existence of a rational point in
  $k$-variety $X$, we need the following lemma, which is probably
  well-known (compare for example the closely related result
  \cite[Theorem~0.3 and Example 0.8]{Milne:ADT}).  We are grateful for J. Bernstein
  for proving this lemma.
\end{subsec}

\begin{lemma}\label{lem:Ext-H^i}
Let  $\Mbul$ be a bounded complex of torsion free finitely generated  (over $\Zz$) discrete $\Gal(\kbar/k)$-modules.
Then for all integers $i$ we have  canonical isomorphisms
$$
\Ext^i(\Mbul,\kbar^\times)=H^i(k,\Hombul_\Zz(\Mbul,\kbar^\times)).
$$
\end{lemma}

Let $\gg$ be a profinite group.
By a $\gg$-module we mean a discrete $\gg$-module.
By a  torsion free finitely generated $\gg$-module we mean a $\gg$-module
which is torsion free and finitely generated over $\Zz$.
Lemma \ref{lem:Ext-H^i} follows from the following Lemma \ref{lem:Ext-Ext}.

\begin{lemma}\label{lem:Ext-Ext}
Let $A$ be a $\gg$-module,  $B$ a $\gg$-module,
and let $\Mbul$ be a complex of torsion free finitely generated  $\gg$-modules.
Then there are canonical isomorphisms
$$
\Ext^i_\gg(A,\Hombul_\Zz(\Mbul,B))=\Ext^i_\gg(A\tensor_\Zz \Mbul,B).
$$
\end{lemma}
To obtain Lemma \ref{lem:Ext-H^i} we just take $A=\Zz,\ B=\kbar^\times$ in Lemma \ref{lem:Ext-Ext}.

\begin{proof}[Proof of Lemma \ref{lem:Ext-Ext}]
First let $M$ be a finitely generated $\gg$-module.
We have a canonical isomorphism
$$
\Hom_\Zz(A,\Hom_\Zz(M,B))=\Hom_\Zz(A\tensor_\Zz M, B).
$$
Taking $\gg$-invariants, we obtain
$$
\Hom_\gg(A,\Hom_\Zz(M,B))=\Hom_\gg(A\tensor_\Zz M, B).
$$

If $\Mbul$ is a complex of torsion free finitely generated $\gg$-modules,
we obtain similarly
$$
\Hombul_\gg(A,\Hombul_\Zz(\Mbul,B))=\Hombul_\gg(A\tensor_\Zz\Mbul, B).
$$

Now let $\Ibul$ be an injective resolution of $B$ in the category of
discrete $\gg$-modules.
Again
$$
\Hombul_\gg(A,\Hombul_\Zz(\Mbul,\Ibul))=\Hombul_\gg(A\tensor_\Zz\Mbul, \Ibul).
$$

By a definition of $\Ext^i$ we have
$$
\sH^i(\Hombul_\gg(A\tensor_\Zz\Mbul, \Ibul))=\Ext^i_\gg(A\tensor_\Zz \Mbul,B).
$$
To prove Lemma \ref{lem:Ext-Ext} it suffices to prove that
$$
\sH^i(\Hombul_\gg(A,\Hombul_\Zz(\Mbul,\Ibul)))=\Ext^i_\gg(A,\Hombul_\Zz(\Mbul,B)).
$$
This follows from the next lemma.
\end{proof}

\begin{lemma}\label{lem:injective-resolution}
$\Hombul_\Zz(\Mbul,\Ibul)$ is an injective resolution of $\Hombul_\Zz(\Mbul,B)$.
\end{lemma}

\begin{proof}
Since $\Mbul$ is a bounded complex of torsion free finitely generated modules,
we see that $\Hombul_\Zz(\Mbul,\Ibul)$ is a resolution of $\Hombul_\Zz(\Mbul,B)$.
This is an injective resolution, since for any
torsion-free finitely generated $\gg$-module $M$
and an injective $\gg$-module $I$, the $\gg$-module
$\Hom_\Zz(M,I)$ is injective (see for example \cite[Lemma~0.5]{Milne:ADT}).
This completes the proofs of Lemmas \ref{lem:injective-resolution},
\ref{lem:Ext-Ext}, and \ref{lem:Ext-H^i}.
\end{proof}


\section{The  extended Picard complex}\label{sec:extended-picard}

\begin{subsec}\label{subsec:triangle}
Let $X$ be a geometrically integral $k$-variety.
Consider the cone
\begin{equation*}
\UPic(\Xbar) := \langle \Gm \to \tau_{\leq 1} R \Gamma_{X/k} \Gm].
\end{equation*}
In more detail: we can represent $\tau_{\leq 1} R \Gamma_{X/k} \Gm$
as a complex in degrees 0 and~1.
We have a homomorphism $i\fcolon \Gm \to H^0(\Xbar, \Gm)$,
which induces a morphism $i_*\fcolon \Gm \to  \tau_{\leq 1} R \Gamma_{X/k} \Gm$.
Then $\UPic(\Xbar)$ is a cone of this map.
Note that the map $i$ is injective, hence $\sH^{-1}(\UPic(\Xbar))=0$,
and $\UPic(\Xbar)[-1]\in\Ob(D^b(k_\etale)^{\geq 1})$.
It follows that $\Hom(\Gm, \UPic(\Xbar)[-1])=0$, so
by \cite[Prop. 1.1.9]{BBD}
$\UPic(\Xbar)$ is defined
up to a canonical isomorphism.
We call $\UPic(\Xbar)$ the \emph{extended Picard complex} of a variety $X$.
We have a canonical distinguished triangle
\begin{equation}\label{eq:triangle-for-UPic}
\Gm\to  \tau_{\leq 1} R \Gamma_{X/k}\Gm  \to  \UPic(\Xbar) \to\Gm[1].
\end{equation}

Note that
\begin{align*}
\sH^0(\UPic(\Xbar)) & = \U(\Xbar), \\
\sH^1(\UPic(\Xbar)) & = \Pic(\Xbar),\\
 \sH^i(\UPic(\Xbar)) & = 0 \text{ for $i \neq 0,1$.}
\end{align*}
Hence $\UPic(\Xbar)$ is indeed a combination of $\Pic(\Xbar)$ and
$U(\Xbar)$.
In particular, if $X$ is projective, then
$\UPic(\Xbar) = \Pic(\Xbar)[-1]$.

The construction of the complex $\UPic(\Xbar)$ is functorial in $X$
in the derived category.
Indeed, a morphism of $k$-varieties $f \fcolon X \to Y$ induces a pull-back
morphism
$f^*\fcolon \tau_{\leq 1} R \Gamma_{Y/k}\Gm
\to \tau_{\leq 1} R \Gamma_{X/k} \Gm$, hence by \cite[Prop. 1.1.9]{BBD}
a canonical morphism
\begin{equation*}
f^* \fcolon \UPic(\Ybar) \to \UPic(\Xbar).
\end{equation*}
\end{subsec}

\begin{subsec}\label{sec:descr-upic}
{\it An explicit presentation of UPic.}
Assume $X$ to be nonsingular.
We write $\Div(\Xbar)$ for the Galois module of divisors on $\Xbar$,
and $\kbar(\Xbar)$ for the rational function field of $\Xbar$.
The divisor map
\begin{equation*}
{\kbar(\Xbar)}^\times \labelto{\divisor} \Div(\Xbar)
\end{equation*}
has kernel equal to $\kbar[\Xbar]^\times$ and cokernel equal to
$\Pic(\Xbar)$.
We write $\KDiv(\Xbar)$ for the complex of Galois modules
$[{\kbar(\Xbar)}^\times \labelto{\divisor} \Div(\Xbar) \rangle$.
We  show below that $\UPic(\Xbar)\simeq \langle \kbar^\times \to \KDiv(\Xbar) ]$

For this, we need the following fact, which should be well-known to
experts, but for which we do not have an explicit reference.
\end{subsec}

\begin{lemma}\label{lem:KDiv-tau}
There is a canonical isomorphism
$$
\KDiv(\Xbar)\isoto\tau_{\leq 1} R \Gamma_{X/k} \Gm\;.
$$
\end{lemma}

To prove Lemma~\ref{lem:KDiv-tau} we need a construction.

\begin{construction}\label{constr:morphism-of-complexes}
Let $K$ be a complex of sheaves on $X$, $K=K^0\to K^1\to\dots$.
\def\GamVk{\Gamma_{X/k}}
We write $\GamVk K=\GamVk K^0\to \GamVk K^1\to\dots$.
By definition of a right derived functor (see for
example~\cite[Def.~III.6.6]{GM}), we have a homomorphism
$$
\GamVk K\to R\GamVk K
$$

Now assume that we have a morphism  $A\to B$  of sheaves on $X$.
Then we have a distinguished triangle
$$
[A\to B\rangle \to A \to B\to [A\to B\rangle[1],
$$
a morphism of triangles

$$
\xymatrix{
[\GamVk A\to \GamVk B\rangle \ar[r] \ar[d]  &\GamVk A \ar[r]\ar[d]
&\GamVk B \ar[r]\ar[d] &[\GamVk A\to \GamVk B\rangle[1] \ar[d] \\
R\GamVk[A\to  B\rangle \ar[r] &R\GamVk A \ar[r] &R\GamVk B \ar[r] &R\GamVk[A\to B\rangle[1]
}
$$
and a commutative diagram with exact rows
$$
\xymatrix@C=.6em{
  \sH^0[\GamVk A\to \GamVk B\rangle \ar[r] \ar[d]&  \GamVk A \ar[r]
  \ar[d]& \GamVk B \ar[r] \ar[d]&
  \sH^1 [\GamVk A\to \GamVk B\rangle[1]\ar[r] \ar[d]& 0 \ar[d]\\
  R^0\GamVk[A\to B\rangle \ar[r] & R^0\GamVk A \ar[r] & R^0\GamVk B
  \ar[r] & R^1\GamVk[A\to B\rangle \ar[r] & R^1\GamVk A
}
$$
\end{construction}

\begin{proof}[Proof of Lemma~\ref{lem:KDiv-tau}]
By \cite[II.1]{Grothendieck:Brauer}
we have a resolution
\begin{equation*}
0 \to \Gm \to \sK^\times_X \to \sDiv_X \to 0
\end{equation*}
of the sheaf $\Gm$ by the sheaf $\sK^\times_X$ of invertible rational
functions
and the sheaf $\sDiv_X$ of Cartier divisors.
Hence we get a canonical isomorphism
\begin{equation}\label{eq:RGammaGm-descr}
R \Gamma_{X/k}\Gm  \simeq R \Gamma_{X/k} [ \sK^\times_X \to \sDiv_X \rangle.
\end{equation}
We have  $R^0 \Gamma_{X/k} \sK^\times_X = {\kbar(\Xbar)}^\times$ and
$R^0 \Gamma_{X/k} \sDiv_X = \Div(\Xbar)$.
Applying Construction~\ref{constr:morphism-of-complexes}
to the morphism of sheaves
$\sK^\times_X \to \sDiv_X$,
we obtain a canonical morphism
\begin{equation}\label{eq:KDiv-RGamma-map}
[ {\kbar(\Xbar)}^\times \labelto{\divisor} \Div(\Xbar)  \rangle \to
R \Gamma_{X/k} \Gm
\end{equation}
and  a commutative diagram with exact rows
\begin{equation*}
\diagram{
\xymatrix@C=.6em{
0 \ar[r] &  \sH^0(\KDiv(\Xbar)) \ar[r] \ar[d] & {\kbar(\Xbar)}^\times \ar[r] \ar[d]^{\simeq}&
\Div(\Xbar) \ar[r] \ar[d]^{\simeq} & \sH^1(\KDiv(\Xbar)) \ar[r] \ar[d] & 0 \\
0 \ar[r] & R^0\Gamma_{X/k} \Gm \ar[r] & R^0 \Gamma_{X/k} \sK^\times_X
\ar[r] & R^0  \Gamma_{X/k} \sDiv_X \ar[r] & R^1\Gamma_{X/k} \Gm \ar[r]
& R^1 \Gamma_{X/k} \sK^\times_X
}
}
\end{equation*}
(we use  the isomorphism \eqref{eq:RGammaGm-descr}).
By Hilbert~90  in Grothendieck's form we have $R^1 \Gamma_{X/k} \sK^\times_X = 0$
(cf.~\cite[II, Lemme 1.6]{Grothendieck:Brauer}).
Hence  the five lemma  gives us that the vertical arrows
$\sH^i(\KDiv(\Xbar))\to R^i\Gamma_{X/k} \Gm $ for $i=0,1$ are isomorphisms.
In other words, the morphism~\eqref{eq:KDiv-RGamma-map}
induces an isomorphism
\begin{equation}\label{eq:KDiv-RGamma-map2}
[ {\kbar(\Xbar)}^\times \labelto{\divisor} \Div(\Xbar)  \rangle \isoto
\tau_{\leq 1} R \Gamma_{X/k} \Gm
\end{equation}
in the derived category.
\end{proof}

\begin{corollary}\label{cor:UPic-KDiv}
There is a canonical isomorphism
$$
\langle \kbar^\times \to \KDiv(\Xbar) ] \isoto \UPic(\Xbar).
$$
\end{corollary}
\begin{proof}
We have a
natural commutative diagram in the derived category
of Galois modules
\begin{equation*}
\xymatrix{
\Gm \ar[d] \ar[r] & \tau_{\leq 1} R \Gamma_{X/k} \Gm \ar[d] \\
\kbar^\times \ar[r] & \KDiv(\Xbar)
}
\end{equation*}
of which the vertical arrows are isomorphisms.
The map $\kbar^\times\to\kbar[\Xbar]^\times=\sH^0(\KDiv(\Xbar))$
is injective.
Now our corollary follows from \cite[Prop. 1.1.9]{BBD}
(similar to the argument in \ref{subsec:triangle}).
\end{proof}

\begin{remark}
Observe that
\begin{equation*}
 \langle \kbar^\times \to \KDiv(\Xbar) ] \simeq [ {\kbar(\Xbar)}^\times/\kbar^\times \to \Div(\Xbar) \rangle.
\end{equation*}
We shall write $\KDiv(\Xbar)/\kbar^\times$ for $[ {\kbar(\Xbar)}^\times/\kbar^\times \to \Div(\Xbar) \rangle$.
Then by Corollary \ref{cor:UPic-KDiv} we have $\KDiv(\Xbar)/\kbar^\times\simeq\UPic(\Xbar)$.
\end{remark}

\begin{remark}
The complex $\KDiv(\Xbar)/\kbar^\times$ is not functorial in $X$ in the category of complexes.
Indeed, neither ${\kbar(\Xbar)}^\times/\kbar^\times$ nor $\Div(\Xbar)$ are functorial in $X$.
\end{remark}

\begin{subsec}\textit{Splitting.}\label{subsec:splitting}

Let $X$ be a nonsingular $k$-variety.
Assume that $X$ has a $k$-point $x$.
We set
\begin{align*}
\Div(\Xbar)_x  &=\{D\in\Div(\Xbar)|\ x\notin\textup{supp}(D)\}\\
{\kbar(\Xbar)}^\times_x  &=\{f\in{\kbar(\Xbar)}^\times|\ \divisor(f)\in \Div(\Xbar)_x \} \\
\KDiv(\Xbar)_x&=[{\kbar(\Xbar)}^\times_x \to \Div(\Xbar)_x\rangle
\end{align*}
By a well-known moving lemma,
the composed map
$$
\Div(\Xbar)_x\to\Div(\Xbar)\to\Pic(\Xbar)
$$
is surjective.
It follows that the morphism of complexes
$$
\KDiv(\Xbar)_x\to\KDiv(\Xbar)
$$
is a quasi-isomorphism.

Set
\begin{gather*}
{\kbar(\Xbar)}^\times_{x,1}=\{f\in{\kbar(\Xbar)}^\times|\ f(x)=1\} \\
\KDiv(X)_{x,1}=[{\kbar(\Xbar)}^\times_{x,1} \to \Div(\Xbar)_x\rangle
\end{gather*}
We have an isomorphism
$$
\kbar^\times\oplus{\kbar(\Xbar)}^\times_{x,1}\isoto {\kbar(\Xbar)}^\times_{x}
$$
given by
$$
(c,f)\mapsto cf\quad \text{ where }\quad c\in\kbar^\times, f\in{\kbar(\Xbar)}^\times_{x,1}.
$$
Hence we obtain an isomorphism
$$
\kbar^\times\oplus\KDiv(\Xbar)_{x,1}\isoto \KDiv(\Xbar)_x\;.
$$
We see that the cone $\langle\kbar^\times\to\KDiv(\Xbar)_x]$
is canonically quasi-isomorphic to $\KDiv(\Xbar)_{x,1}$.
Thus $\UPic(\Xbar)\simeq \KDiv(\Xbar)_{x,1}$.

Let $f\fcolon X\to Y$ be a morphism of nonsingular $k$-varieties,
and let $x\in X(k)$.  Set $y=f(x)\in Y(k)$.  Then we have a
morphism of complexes
$$
f^*\fcolon\KDiv(\Ybar)_{y,1}\to\KDiv(\Xbar)_{x,1}\;.
$$
We see that the complex $\KDiv(\Xbar)_{x,1}$ is functorial in $(X,x)$ in the category of complexes.

\begin{lemma}\label{lem:splitting}
Let $X$ be a nonsingular $k$-variety having a $k$-point $x$.
Then the triangle~\eqref{eq:triangle-for-UPic} of \ref{subsec:triangle} splits,
i.e.  the third morphism $\UPic(\Xbar)\to\Gm[1]$
in this triangle is 0.
\end{lemma}
\begin{proof}
If $X$ has a $k$-point $x$, then the triangle \eqref{eq:triangle-for-UPic}
is isomorphic to the {split} triangle
$$
\kbar^\times\to \kbar^\times\oplus\KDiv(\Xbar)_{x,1}\to\KDiv(\Xbar)_{x,1}\to\kbar^\times[1]
$$
with obvious morphisms, where the third morphism is 0.
Hence  the third morphism  in the triangle \eqref{eq:triangle-for-UPic}  is~0.
\end{proof}
\end{subsec}

The lemma shows that the triangle~\eqref{eq:triangle-for-UPic}
can provide a cohomological obstruction to the existence of a
$k$-rational point.

\begin{definition}\label{def:elementary-obstruction}
Let $X$ be a nonsingular variety over $k$.
We define the \emph{elementary obstruction}
\[ e(X) \in \Ext^1(\UPic(\Xbar), \Gm)\]
to be the class $e(X)$ of the
triangle~\eqref{eq:triangle-for-UPic}.
\end{definition}

\begin{subsec}
We call $e(X)$ the elementary obstruction, because it is  closely related
to the {elementary   obstruction}
$\ob(X)\in \Ext^1({\kbar(\Xbar)}^\times/\kbar^\times,\kbar^\times)$
of Colliot-Th\'el\`ene and Sansuc \cite[D\'ef.~2.2.1]{CT-Sansuc}.
Indeed, by definition $\ob(X)$ is the class of the extension
\[ 0 \to \kbar^\times \to {\kbar(\Xbar)}^\times \to
{\kbar(\Xbar)}^\times/\kbar^\times \to 0,
\]
whereas under the identification $\UPic(\Xbar) \simeq
\KDiv(\Xbar)/\kbar^\times$ of Corollary~\ref{cor:UPic-KDiv}, $e(X)$ is the
extension class of the triangle associated to the short exact sequence
of complexes
\[ 0 \to \kbar^\times \to \KDiv(\Xbar) \to \KDiv(\Xbar)/\kbar^\times \to 0. \]
Hence $e(X)$ is the image of the class $\ob(X)$ under the homomorphism
\begin{equation}\label{eq:ob-e-map}
 \Ext^1({\kbar(\Xbar)}^\times/\kbar^\times, \kbar^\times) \to
\Ext^1(\KDiv(\Xbar)/\kbar^\times, \kbar^\times)
\end{equation}
induced by the
natural map $\KDiv(\Xbar)/\kbar^\times \to
{\kbar(\Xbar)}^\times/\kbar^\times$.
\end{subsec}

\begin{lemma}
For a nonsingular $k$-variety $X$,
we have $e(X)=0$ if and only if $\ob(X)=0$.
\end{lemma}
\begin{proof}
The homomorphism~\eqref{eq:ob-e-map}
fits into an exact sequence of $\Ext$-groups
\begin{equation}\label{eq:exact-Ext}
  \Ext^1(\Div(\Xbar),\kbar^\times) \to
  \Ext^1({\kbar(\Xbar)}^\times/\kbar^\times, \kbar^\times) \to
  \Ext^1\left([ {\kbar(\Xbar)}^\times/\kbar^\times \to
  \Div(\Xbar) \rangle, \kbar^\times\right).
\end{equation}
induced by the exact sequence of complexes
$$
0\to \Div(\Xbar)[-1]
\to [{\kbar(\Xbar)}^\times/\kbar^\times\to\Div(\Xbar)\rangle
\to {\kbar(\Xbar)}^\times/\kbar^\times
\to 0
$$

Since $\Div(\Xbar)$ is a direct sum of permutation modules,
Lemma \ref{lem:Ext-H^i} gives that $\Ext^1(\Div(\Xbar),\kbar^\times)$
is a direct product of the $H^1$-groups of quasi-trivial tori,
hence $\Ext^1(\Div(\Xbar),\kbar^\times)=0$,
so we see from the exact sequence~\eqref{eq:exact-Ext}
that the homomorphism~\eqref{eq:ob-e-map} is injective, from which the
statement follows.
\end{proof}


Now we investigate how UPic changes under open embeddings.

\begin{proposition}\label{prop:zariski-open-upic}
Let $X \subset Y$ be an open $k$-subvariety of a nonsingular $k$-variety $Y$.
Let $j\fcolon X\into Y$ denote the inclusion map.
Then we have a distinguished triangle
\begin{equation*}
  \UPic(\Ybar) \labelto{j^*} \UPic(\Xbar) \to \cycles^1_{Y - X}\to \UPic(\Ybar)[1]
\end{equation*}
where $\cycles^1_{Y - X}$ is the permutation module of
divisors in the complement of $\Xbar$ in $\Ybar$.
\end{proposition}

\begin{proof}
Clearly we have a short exact sequence of complexes
$$
0\to  \cycles^1_{Y-X}[-1] \to  \KDiv(\Ybar)/\kbar^\times \labelto{j^*}
\KDiv(\Xbar)/\kbar^\times\to 0,
$$
whence we obtain distinguished triangles
$$
  \cycles^1_{Y-X}[-1] \to  \KDiv(\Ybar)/\kbar^\times \labelto{j^*}
\KDiv(\Xbar)/\kbar^\times\to \cycles^1_{Y-X}
$$
and
$$
\KDiv(\Ybar)/\kbar^\times \labelto{j^*} \KDiv(\Xbar)/\kbar^\times
\labelto{}
\cycles^1_{Y-X}  \labelto{} (\KDiv(\Ybar)/\kbar^\times)[1]
$$
\end{proof}

\begin{remark}
Let $X\subset Y$ be an open $k$-subvariety of a nonsingular
complete $k$-variety $Y$.
Proposition \ref{prop:zariski-open-upic} implies that $\UPic(\Xbar)$ is
non-canonically isomorphic to the fibre $\left[\cycles^1_{Y-X}\to\Pic(\Ybar)\right\rangle$.
Skorobogatov actually gave a canonical isomorphism in the derived category
$\UPic(\Xbar)\isoto \left[\cycles^1_{Y-X}\to\Pic(\Ybar)\right\rangle$
(cf. \cite[Rem. B.2.1(2)]{CT06}).
\end{remark}

By $\Sha_\omega^i(k, M)$ we denote the subgroup of $H^i(k, M)$ of
elements that map to zero in $H^i(\gamma, M)$ for every closed
procyclic subgroup $\gamma \subset \Gal(\kbar/k)$.  Recall that for
a permutation module $P$ we have  $H^1(k, P) = 0$ and
$\Sha_\omega^2(k, P) = 0$ (cf.\
\cite[1.2.1]{Borovoi-Kunyavskii:formulas}).

\begin{corollary}\label{cor:h2-upic-open-dense}
Let $X \subset Y$ be an open $k$-subvariety of a nonsingular $k$-variety $Y$.
Then the restriction map
$\UPic(\Ybar) \to \UPic(\Xbar)$
induces an injection
\begin{equation*}
H^2(k, \UPic(\Ybar)) \into H^2(k, \UPic(\Xbar))
\end{equation*}
and an isomorphism
\begin{equation*}
\Sha_\omega^2(k, \UPic(\Ybar)) \isoto \Sha_\omega^2(k, \UPic(\Xbar)).
\end{equation*}
\end{corollary}
\begin{proof}
By  Proposition~\ref{prop:zariski-open-upic} we have an exact sequence
$$
H^1(k,\cycles^1_{Y-X})\to H^2(k, \UPic(\Ybar)) \to H^2(k, \UPic(\Xbar))\to H^2(k,\cycles^1_{Y-X}),
$$
where $\cycles^1_{Y-X}$ is a permutation Galois module.
Now the injectivity of the two maps follows from
the vanishing of $H^1(k, \cycles^1_{Y - X})$.
The surjectivity of the $\Sha^2_\omega$-map follows from the vanishing
of $\Sha^2_\omega(k, \cycles^1_{Y - X})$
and an easy diagram chase.
\end{proof}

\begin{corollary}\label{cor:h1-pic-open-dense}
Let $X\subset Y$ be an open $k$-subvariety of a nonsingular
complete $k$-variety $Y$.
Then the restriction map
$\Pic(\Ybar)[-1] = \UPic(\Ybar) \to \UPic(\Xbar)$
induces an injection
\begin{equation*}
H^1(k, \Pic(\Ybar)) \into H^2(k, \UPic(\Xbar))
\end{equation*}
and an isomorphism
\begin{equation*}
\Sha_\omega^1(k, \Pic(\Ybar)) \isoto \Sha_\omega^2(k, \UPic(\Xbar)).
\end{equation*}%
\end{corollary}

\begin{corollary}\label{cor:restr-elementary}
Let $X \subset Y$ be
an open $k$-subvariety of a nonsingular $k$-variety $Y$.
Let $j\fcolon X\into Y$ denote the inclusion map.
Then the induced map
\begin{equation*}
j_* \fcolon \Ext^1(\UPic(\Xbar),\Gm)
                   \to \Ext^1(\UPic(\Ybar),\Gm)
\end{equation*}
is injective.
In particular, the elementary obstruction  $e(X)$
vanishes if and only if  $e(Y)$ vanishes.
\end{corollary}
\begin{proof}
  Applying  the functor $\Ext$
  to the distinguished triangle of
  Proposition~\ref{prop:zariski-open-upic},
  we obtain an exact sequence
  \[ \Ext^1(\cycles^1_{Y - X},\Gm) \to
     \Ext^1(\UPic(\Xbar),\Gm) \labelto{j_*}
     \Ext^1(\UPic(\Ybar),\Gm).
  \]
By Lemma \ref{lem:Ext-H^i} $\Ext^1(\cycles^1_{Y - X},\Gm)=H^1(k,P)$,
where $P$ is the $k$-torus such that $\Chi^*(P)=\cycles^1_{Y - X}$.
Since $\cycles^1_{Y - X}$ is a permutation module,
we see that $P$ is a quasi-trivial torus,
hence $H^1(k,P)=0$, and therefore the homomorphism $j_*$ is injective.
\end{proof}

\begin{subsec}
{\it UPic, the Picard group and the Brauer group.}
Let $X$ be a nonsingular variety over $k$.
Let $\Br(X) = H^2(X, \Gm)$ denote the (cohomological)
Brauer group of $X$, let
$\Br_1(X)$ denote the kernel of the map $\Br(X) \to \Br(\Xbar)$,
and let $\Bra(X)$ denote the cokernel of the map $\Br(k)\to\Br_1(X)$.

We have equalities
\begin{align*}
\Pic(X) & = H^1(X, \Gm) = H^1(k, R \Gamma_{X/k} \Gm) = H^1(k, \tau_{\leq 1} R
\Gamma_{X/k} \Gm) \\
\Br(X) & = H^2(X, \Gm) = H^2(k, R \Gamma_{X/k} \Gm) =
H^2(k, \tau_{\leq 2} R \Gamma_{X/k} \Gm).
\end{align*}
{}From the distinguished triangle
$$
\tau_{\leq 1} R \Gamma_{X/k} \Gm \to \tau_{\leq 2} R \Gamma_{X/k} \Gm
\to  R^2 \Gamma_{X/k} \Gm[-2] \to \tau_{\leq 1} R \Gamma_{X/k} \Gm[1]
$$
we obtain a Galois cohomology exact sequence
$$
0\to H^2(k, \tau_{\leq 1} R \Gamma_{X/k} \Gm)\to H^2(k, \tau_{\leq 2} R \Gamma_{X/k} \Gm)
\to H^0(k,  R^2 \Gamma_{X/k} \Gm).
$$
Since $H^2(k, \tau_{\leq 2} R \Gamma_{X/k} \Gm) = \Br(X)$, and
$H^0(k, R^2 \Gamma_{X/k} \Gm) = \Br(\Xbar)^{\Gal(\kbar/k)}$,
it follows that
\begin{equation*}
\Br_1(X) =  H^2(k, \tau_{\leq 1} R \Gamma_{X/k} \Gm ).
\end{equation*}
\end{subsec}
\begin{proposition}\label{prop:pic-brauer-upic-sequence}
Let $X$ be a nonsingular variety over $k$.
We have an exact sequence
\begin{multline*}
0 \to \Pic(X) \to H^1(k, \UPic(\Xbar)) \to \Br(k) \to \Br_1(X)
\to H^2(k, \UPic(\Xbar)) \to H^3(k, \Gm),
\end{multline*}
in which the homomorphisms $H^1(k, \UPic(\Xbar)) \to \Br(k)$ and
$ H^2(k, \UPic(\Xbar)) \to H^3(k, \Gm)$ are zero when $X(k) \neq
\emptyset$.
\end{proposition}
\begin{proof}
We obtain the exact sequence by taking Galois cohomology of the
triangle~\eqref{eq:triangle-for-UPic} of \ref{subsec:triangle} and applying
Hilbert's Theorem~90 to the term $H^1(k, \Gm)$.
For the case $X(k) \neq \emptyset$ we apply Lemma~\ref{lem:splitting}.
\end{proof}

\begin{corollary}\label{cor:pic-brauer-upic}
Let $X$ be a smooth geometrically integral variety over $k$.
\begin{theoremlist}
\item\label{corpart:pic-upic}
There is a canonical injection
$$
\Pic(X)\into H^1(k,\UPic(\Xbar))
$$
which is an isomorphism if $X(k) \neq \emptyset$ or if $\Br(k) =0$.
\item\label{corpart:brauer-upic}
There is a canonical injection
\begin{equation*}
  \Bra(X) \into H^2(k,\UPic(\Xbar) )
\end{equation*}
which is an isomorphism if $X(k) \neq \emptyset$ or $H^3(k,\Gm) = 0$.
\end{theoremlist}
\end{corollary}


\section{Picard groups, invertible functions, and the algebraic
  fundamental group}

\begin{subsec}
Let $G$ be a connected linear algebraic $k$-group.
 As in \cite{Borovoi:Memoir} we write
$G\unipotent \subset G$ for the unipotent radical of $G$,
$G\reductive$ for the reductive group $G/G\unipotent$,
$G\semisimple$ for the derived group  of $G\red$ (it is semisimple),
$G\tor$ for the torus $G\reductive/G\semisimple$, and
$G\sc$ for the universal covering of $G\semisimple$ (it is simply connected).
The composed map
\begin{equation*}
\rho \fcolon G\sc \onto G\semisimple \into G\reductive
\end{equation*}
has finite kernel
\[
Z := \ker \rho,
\]
which is central in $G\sc$,
and the cokernel of $\rho$ is equal
 to the torus $G\tor$.
We write
$$
\Chi^*(\Gbar)=\Hom_{\kbar}(G,\Gm)
$$
for the character group of $G$.
We have
$$
\Chi^*(\Gbar)=\Chi^*(\Gbar\tor).
$$
For a torus $T$ we write
$$
\Chi_*(\Tbar)=\Hom_{\kbar}(\Gm,T)
$$
for the cocharacter group of $T$.
Note that the underlying abelian groups of the Galois modules
$\Chi^*(\Gbar)$ and $\Chi_*(\Tbar)$ are free.

As in \cite{Borovoi:Memoir} we define the algebraic fundamental group
$\pi_1(\Gbar)$ as follows.
Let $T\subset G\red$ be a maximal torus.
Set $T\sc=\rho^{-1}(T)$, it is a maximal torus in $G\sc$.
Set
$$
\pi_1(\Gbar)=\Chi_*(\Tbar)/\rho_*\Chi_*(\Tbar\sc).
$$
It is a Galois module; it does not depend on the choice of $T\subset G$;
it does not change under inner twistings of $G$.
It follows from the definition, that $\pi_1(\Gbar)=\pi_1(\Gbar\red)$.

We define the derived dual to $\pi_1(\Gbar)$ by
$$
\pi_1(\Gbar)^D = [ \Chi^*(\Tbar) \labelto{\rho^*} \Chi^*(\Tbar\sc)\rangle.
$$
In this section we shall recall classical results that give isomorphisms
\begin{align*}
\sH^0(\UPic(\Gbar)) = U(\Gbar) & = \ker \rho^* = \sH^0(\pi_1(\Gbar)^D), \\
\intertext{and}
\sH^1(\UPic(\Gbar)) = \Pic(\Gbar) & = \coker \rho^* = \sH^1(\pi_1(\Gbar)^D). \\
\end{align*}

\end{subsec}

\begin{lemma}[Rosenlicht]\label{lem:rosenlicht}
For a connected linear algebraic group $G$ over a perfect field $k$,
the obvious map $\Chi^*(\Gbar)\to \U(\Gbar)$ is an isomorphism which is functorial in $G$.
\end{lemma}
\begin{proof}
See \cite{Rosenlicht}, or \cite[Cor. 2.2]{Fossum-Iversen}, or \cite[Prop. 1.2]{KKV}
\end{proof}

\begin{corollary}
  For a connected linear  $k$-group $G$ we have a canonical isomorphism
\[ U(\Gbar) \iso \sH^0(\pi_1(\Gbar)^D). \]
\end{corollary}
\begin{proof}
Clearly, $\Chi^*(\Gbar\tor) \simeq \Chi^*(\Gbar)$, hence
$\Chi^*(\Gbar\tor) \simeq U(\Gbar)$.
On the other hand, the identification $G\tor = T/\rho(T\sc)$
gives an isomorphism
\[ \Chi^*(\Gbar\tor) \simeq
\ker[\Chi^*(\Tbar) \to \Chi^*(\Tbar\sc)]
= \sH^0(\pi_1(\Gbar)^D). \]
\end{proof}

We shall now consider the identification
$\sH^1(\UPic(\Gbar)) = \sH^1(\pi_1(\Gbar)^D).$
We first make a reduction to $G\sss$ using
the following lemma of Fossum--Iversen and Sansuc.

\begin{lemma}\label{lem:Fossum-Iversen-Sansuc}
Let $1\to G'\to G\to G''\to 1$ be an exact sequence of connected
linear $k$-groups.
Then we have an exact sequence
$$
0\to\Chi^*(\Gbar'')\to\Chi^*(\Gbar)\to\Chi^*(\Gbar')\to \Pic(\Gbar'')\to\Pic(\Gbar)\to\Pic(\Gbar')\to 0.
$$
\end{lemma}
\begin{proof}
See \cite[(6.11.4)]{Sansuc:brauer-gal}.
In the case when $H^1(K, G')=0$ for any extension $K$ of $\kbar$,
this exact sequence was obtained in \cite[Prop. 3.1]{Fossum-Iversen}.
\end{proof}

\begin{corollary}\label{cor:pic-mod-uni-ss}
Let $G$ be a linear algebraic group over $k$.
Then the canonical maps $r\fcolon G \to G\reductive$
and $G\sss\to G\reductive$
induce a natural isomorphism
\begin{equation*}
\Pic(\Gbar\sss) \iso \Pic(\Gbar).
\end{equation*}
\end{corollary}
\begin{proof}
We first apply
Lemma~\ref{lem:Fossum-Iversen-Sansuc}
to the short exact sequence
\[ 1 \to G\uu \to G \to G\red \to 1, \]
and then to the short exact sequence
\[ 1 \to G\sss \to G\red \to G\tor \to 1. \]
using the fact that $\Chi^*(\Gbar\uu) = 0$,  $\Pic(\Gbar\uu)=0$,
and $\Pic(G\tor) = 0$.
\end{proof}

We need the following construction of
\cite[\S2]{Popov} (see also  \cite[p.~275]{Fossum-Iversen} and \cite[Example 2.1]{KKLV89}).

\begin{construction}\label{constr:FI}
Let $G$ be a connected linear $k$-group.
Let $H\subset G$ be a $k$-subgroup, not necessarily connected.
Set $X=G/H$.
We construct a morphism of Galois modules
$$
c\fcolon \Chi(\Hbar)\to\Pic(\Xbar)
$$
as follows.
Let $\chi\in\Chi(\Hbar)$.
Consider the embedding
$$
\Hbar\into\Gbar\times\Gmbar,\quad h\mapsto (h,\chi(h)^{-1}).
$$
Set $\Ybar=(\Gbar\times\Gmbar)/\Hbar$,
this quotient exists by Chevalley's theorem, see for example \cite[Thm. 5.5.5]{Springer}.
We have a canonical map $\Ybar\to\Xbar=\Gbar/\Hbar$.
Clearly $\Ybar$ is a torsor under $\Gmbar$ over $\Xbar$,
which admits a local section (in Zariski topology) by Hilbert 90.
Since  the group $\Gbar\times\Gmbar$ acts transitively on $\Ybar$ and $\Xbar$,
we conclude that the torsor $\Ybar\to\Xbar$ is locally trivial in Zariski topology.
From the principal $\Gmbar$-bundle $\Ybar$ we construct
(using the transition functions of $\Ybar$)
a \emph{linear} bundle on $\Xbar$ which we denote by $L(\chi)$.
Alternatively, we can construct $L(\chi)$ directly
as the quotient $(\Gbar\times\Gak)/\Hbar$ of $\Gbar\times \Gak$
under the right action of $\Hbar$
given by
$$
(g,a)\cdot h=(gh, a \chi(h)^{-1}),
$$
where $g\in \Gbar,\ a\in\Gak,\ h\in \Hbar$.

We denote by $c(\chi)$ the class of $L(\chi)$ in $\Pic(\Xbar)$.
In terms of divisor classes, this means the following.
Let $\psi_X$ be a rational section of $L(\chi)$.
Set $D=\divisor(\psi_X)$.
We set $c(\chi)=\cl(D)\in\Pic(\Xbar)$,
where $\cl(D)$ denotes the class of the divisor $D$.

Note that a rational section $\psi_X$ of $L(\chi)$ over $\Xbar$
lifts canonically to a rational function $\psi_G$ on $\Gbar$.
Namely, the graph of $\psi_G$ in $\Gbar\times\Gak$ is the preimage
of the graph of $\psi_X$ in $L(\chi)$ with respect to the quotient map
$\Gbar\times\Gak\to L(\chi)$.
\end{construction}

\begin{lemma}\label{lem:KKV} \cite[Thm. 4]{Popov},  \cite[Prop. 3.2]{KKV}
Let $G$ be a connected linear $k$-group, and let $H$ be a $k$-subgroup of $G$
(not necessarily connected).
Then the sequence
$$
\Chi^*(\Gbar)\to\Chi^*(\Hbar)\labelto{c}\Pic(\Gbar/\Hbar)\to\Pic(\Gbar)
$$
is exact.
\end{lemma}

\begin{corollary}\label{cor:FI}
Let $G$ be a connected semisimple $k$-group.
We regard $G$ as a homogeneous space $G=G\sc/Z$, where $Z=\ker\rho$.
Then we have an isomorphism
$$
c\fcolon \Chi^*(\Zbar)\isoto\Pic(\Gbar),
$$
where $c$ is the homomorphism of Construction \ref{constr:FI}.
\end{corollary}
\begin{proof}
The corollary follows from Lemma \ref{lem:KKV}.
We use the facts that $\Chi^*(\Gbar\sc)=0$ and $\Pic(\Gbar\sc)=0$.
\end{proof}

\begin{remark}
The equality $\Pic(\Gbar\sc)=0$ and
the existence of an isomorphism $\Chi^*(\Zbar)\simeq\Pic(\Gbar)$ for a semisimple $k$-group $G$
were proved by Voskresenski\u\i\  \cite{Vos},  Fossum and Iversen \cite[Cor. 4.6]{Fossum-Iversen},
and Popov \cite{Popov} (see also \cite[4.3]{Vos98}).
\end{remark}

\begin{corollary}
\label{cor:pic-g-chi-z}
For any connected linear $k$-group $G$
we have a canonical isomorphism
\[ \Chi^*(\Zbar) \isoto \Pic(\Gbar) \]
where $Z=\ker\rho$.
\end{corollary}
\begin{proof}
By Corollary~\ref{cor:pic-mod-uni-ss}, we have an isomorphism
$\Pic(\Gbar\sss) \isoto \Pic(\Gbar)$.
By Corollary \ref{cor:FI} we have an isomorphism
$\Chi^*(\Zbar) \isoto \Pic(\Gbar\sss)$.
\end{proof}

\begin{corollary}
For any connected linear $k$-group $G$
we have a canonical isomorphism
$\Pic(\Gbar)\simeq\sH^1(\pi_1(\Gbar)^D)$.
\end{corollary}
\begin{proof}
Indeed,
$$
\Chi^*(\Zbar)= \Chi^*(\ker \rho) = \sH^1([\Chi^*(T) \labelto{\rho^*}
\Chi^*(T\sc) \rangle = \sH^1(\pi_1(\Gbar)^D).
$$
\end{proof}

\section{The Comparison Theorem}

We shall now construct a canonical isomorphism
$$
\varkappa_G\colon \UPic(\Gbar)\to\pi_1(\Gbar)^D
$$
for any connected linear $k$-group $G$.

We first make a reduction to the reductive case.

\begin{lemma}\label{lem:upic-mod-uni}
Let $G$ be a connected linear algebraic group over $k$.
Then the canonical homomorphism $r\colon G \to G\reductive$
induces an isomorphism
\begin{equation*}
r^*\colon \UPic(\Gbar\reductive) \isoto \UPic(\Gbar).
\end{equation*}
\end{lemma}
\begin{proof}
By Lemma \ref{lem:Fossum-Iversen-Sansuc} we have an exact sequence
$$
0\to \Chi^*(\Gbar\red)\to\Chi^*(\Gbar)\to\Chi^*(\Gbar\uu)\to \Pic(\Gbar\red)
\to\Pic(\Gbar)\to\Pic(\Gbar\uu)\to 0
$$
where $\Chi^*(\Gbar\uu)=0$ and $\Pic(\Gbar\uu)=0$.
It follows that the map $r$ induces isomorphisms
$\Chi^*(\Gbar\reductive) \isoto \Chi^*(\Gbar)$ and
$\Pic(\Gbar\reductive) \isoto \Pic(\Gbar)$.
We see that the  morphism
$r^*\colon \UPic(\Gbar\reductive) \to \UPic(\Gbar)$
induces isomorphisms on $\sH^0$ and $\sH^1$, hence it is an isomorphism.
\end{proof}

\begin{lemma}\label{lem:torus-upic}
  For any torus $T$ over  $k$ we have
  $\UPic(\Tbar) \simeq \Chi^*(\Tbar)$.
\end{lemma}
\begin{proof}
This follows from Rosenlicht's lemma (Lemma \ref{lem:rosenlicht}), since
$\Pic(\Tbar) = 0$.
\end{proof}

\begin{lemma}\label{lem:upic-ss-sc}
Let $G$ be a connected linear algebraic group over $k$ such that
$G\semisimple$ is simply connected, then we have canonical isomorphisms
$\UPic(\Gbar) = \Chi^*(\Gbar) = \Chi^*(\Gbar\tor)$.
In particular, $\UPic(\Gbar) = 0$ if $G$ is semi-simple and simply connected.
\end{lemma}

\begin{proof}
We have   $\Chi^*(\Gbar) = \Chi^*(\Gbar\tor)$ (for any connected $G$).
By Lemma~\ref{cor:pic-g-chi-z} $\Pic(\Gbar)=0$,
hence $\UPic(\Gbar)=\Chi^*(\Gbar)$.
\end{proof}

\begin{subsec}\label{subsec:morphisms-of-complexes}
In this subsection  and  the next one we identify $\UPic(\Gbar)$ with $\KDiv(\Gbar)_{e,1}$
as in \ref{subsec:splitting},
where $e$ is the unit element of $G$.
We write $\KDiv(\Gbar)_{1}$ for $\KDiv(\Gbar)_{e,1}$.
Note that $G\mapsto \KDiv(\Gbar)_{1}$ is a functor
from the category of connected linear $k$-groups
to the category of complexes of Galois modules.

For a maximal torus $T$ in a connected reductive $k$-group $G$
we have a commutative diagram
$$
\xymatrix{
T\sc \ar[r]^{\rho_T} \ar[d]_{i\sc}  &T \ar[d]^i  \\
G\sc \ar[r]^{\rho}    &G
}
$$
(where $i$ is the inclusion homomorphism), hence a commutative diagram of complexes
\begin{equation}\label{eq:CD-GT}
\xymatrix{
\KDiv(\Gbar)_1 \ar[r]^{\rho^*} \ar[d]_{i^*} &\KDiv(\Gbar\sc)_1 \ar[d]^{(i\sc)^*} \\
\KDiv(\Tbar)_1 \ar[r]^{\rho_T^*}            &\KDiv(\Tbar\sc)_1
}
\end{equation}
and a morphism of complexes
$$
\lambda=i^*\oplus(i\sc)^*\colon[\KDiv(\Gbar)_1\to\KDiv(\Gbar\sc)_1\rangle\to [\KDiv(\Tbar)_1\to\KDiv(\Tbar\sc)_1\rangle.
$$

Consider the fibre $[\KDiv(\Gbar)_1\to\KDiv(\Gbar\sc)_1\rangle$.
The canonical morphism
$$
[\KDiv(\Gbar)_1\to\KDiv(\Gbar\sc)_1\rangle\to\KDiv(\Gbar)_1
$$
is an isomorphism in the derived category,
because $\KDiv(\Gbar\sc)_1\simeq 0$ by Lemma~\ref{lem:upic-ss-sc}.

Consider the fibre $[\KDiv(\Tbar)_1\to\KDiv(\Tbar\sc)_1\rangle$.
The commutative diagram of complexes
\begin{equation*}
\xymatrix{
{\Chi^*(\Tbar)}\ar[r]\ar[d] &\Chi^*(\Tbar\sc)\ar[d]\\
\KDiv(\Tbar)_1\ar[r]          &\KDiv(\Tbar\sc)_1
}
\end{equation*}
in which the vertical arrows are isomorphisms in the derived category,
induces a morphism of complexes
$$
[\Chi^*(\Tbar)\to\Chi^*(\Tbar\sc)\rangle\to
[\KDiv(\Tbar)_1\to\KDiv(\Tbar\sc)_1\rangle
$$
which is an isomorphism in the derived category.
\end{subsec}

\begin{construction}\label{constr-morphism}
For a reductive $k$-group $G$
we define a morphism
$\varkappa_G\colon \UPic(\Gbar)\to\pi_1(\Gbar)^D$
as the  composition
\begin{align*}
&\UPic(\Gbar)\isoto[\KDiv(\Gbar)_1\to\KDiv(\Gbar\sc)_1\rangle\labelto{\lambda}
[\KDiv(\Tbar)_1\to\KDiv(\Tbar\sc)_1\rangle \\
&\isoto [\Chi^*(\Tbar)\to\Chi^*(\Tbar\sc)\rangle=\pi_1(\Gbar)^D.
\end{align*}
For a general connected linear $k$-group $G$ (not necessarily reductive) we define
$\varkappa_G\colon \UPic(\Gbar)\to\pi_1(\Gbar)^D$
as the composition
$$
\UPic(\Gbar)\isoto\UPic(\Gbar\red)\to\pi_1(\Gbar\red)^D=\pi_1(\Gbar)^D.
$$
\end{construction}

\begin{remark}
Using \cite[Prop. 1.1.9]{BBD}, one can show that all the fibres and morphisms of fibres
that we  defined in
\ref{subsec:morphisms-of-complexes} and \ref{constr-morphism}
using an explicit representation by complexes and morphisms of complexes,
do not depend on this representation.
We use the fact that $\Pic(\Gbar)$ is a torsion group (by Corollary \ref{cor:pic-g-chi-z}),
hence
$$
\Hom_{D^b(k_\etale)}(\UPic(\Gbar)[1], \UPic(\Tbar\sc))=\Hom(\Pic(\Gbar), \Chi^*(\Tbar\sc))=0,
$$
and so \cite[Prop. 1.1.9]{BBD} can be applied.
\end{remark}

\begin{subsec}\label{subsec:functoriality} \emph{Functoriality.}
Let $\varphi\colon G_1\to G_2$ be a homomorphism of connected linear $k$-groups.
Consider the induced homomorphisms $\varphi\red\colon G_1\red\to G_2\red$
and $\varphi\sc\colon G_1\sc\to G_2\sc$.
Choose maximal tori $T_1\subset G_1\red$ and $T_2\subset G_2\red$
such that $\varphi\red(T_1)\subset T_2$.
Let $T_1\sc$ (resp. $T_2\sc$) be the preimage of $T_1$ in $G_1\sc$ (resp. of $T_2$ in $G_2\sc$).
We have homomorphisms $\varphi_*\colon T_1\to T_2$ and $\varphi_*\sc\colon T_1\sc\to T_2\sc$.
We obtain a commutative diagram in the derived category
\def\XX{{\Chi^*}}
$$
\xymatrix{
{\KDiv(\Gbar_2)_1}\ar[r]^{\varphi^*}\ar[d]_{\varkappa_{G_2}}  &\KDiv(\Gbar_1)_1\ar[d]_{\varkappa_{G_1}}  \\
{[\XX(\Tbar_2)\to\XX(\Tbar_2\sc)\rangle}\ar[r]^{\varphi^*}  &[\XX(\Tbar_1)\to\XX(\Tbar_1\sc)\rangle
}
$$
Thus the morphism $\varkappa_G\colon\UPic(\Gbar)\to\pi_1(\Gbar)^D$ is functorial in $G$.
\end{subsec}

The following theorem is the main result of this paper.

\begin{theorem}\label{thm:main}
With notation as above, for a connected linear algebraic group $G$
over a field $k$ of characteristic 0,
the canonical morphism $\varkappa_G\fcolon\UPic(\Gbar)\to\pi_1(\Gbar)^D$
is an isomorphism.
\end{theorem}

Before proving  the theorem, let us first mention two corollaries.

\begin{corollary}\label{cor:upic-dual-tori}
The canonical isomorphism $\varkappa_G$ induces  canonical
isomorphisms
$$ \Ext^i(\UPic(\Gbar), \Gm) \simeq H^i_\ab(k, G),$$
where $H^i_\ab(k,G):=H^i(k,\langle T\sc \to T ])$.
\end{corollary}
\begin{proof}
  By Theorem \ref{thm:main} $\Ext^i(\UPic(\Gbar), \Gm)=\Ext^i(\pi_1(\Gbar)^D, \Gm)$.
By Lemma \ref{lem:Ext-H^i} $\Ext^i(\pi_1(\Gbar)^D, \Gm)=H^i(k,\langle T\sc \to T ])$.
\end{proof}

\begin{corollary}
Let $1\to G'\to G\to G''\to 1$ be an exact sequence of connected
linear $k$-groups.
Then we have a distinguished triangle
$$
\UPic(\Gbar'')\to\UPic(\Gbar)\to\UPic(\Gbar')\to\UPic(\Gbar'')[1].
$$
\end{corollary}

\begin{proof}
Indeed, by \cite[Lemma 1.5]{Borovoi:Memoir}
(see also \cite[Lemma 3.7]{Borovoi-Kunyavskii:arithmetical})
we have an exact sequence
$$
0\to \pi_1(\Gbar')\to\pi_1(\Gbar)\to\pi_1(\Gbar'')\to 0,
$$
hence a distinguished triangle
$$
\pi_1(\Gbar'')^D\to\pi_1(\Gbar)^D\to\pi_1(\Gbar')^D\to\pi_1(\Gbar'')^D[1],
$$
and the assertion follows from Theorem \ref{thm:main}.
\end{proof}

\begin{remark}
This triangle strengthens Lemma~\ref{lem:Fossum-Iversen-Sansuc}.
Note that it does not give a new proof, since the lemma
was used in the proof of  Theorem~\ref{thm:main},
hence in the proof of
this corollary.
\end{remark}

\def\KD{\text{KDiv}}

\begin{subsec}
{\it Proof of Theorem \ref{thm:main}.}
We may and shall assume that $G$ is reductive.
Recall that
$$
\KDiv(\Gbar)_1=[\kbar(\Gbar)^\times_{e,1}\to\Div(\Gbar)_e\rangle,
$$
where
\begin{gather*}
\kbar(\Gbar)^\times_{e,1}=\{f\in\kbar(\Gbar)^\times|\ f(e)=1\},\\
\Div(\Gbar)_e=\{D\in\Div(\Gbar)|\ e\notin\textup{supp}(D)\}.
\end{gather*}
From the diagram \eqref{eq:CD-GT}
we obtain a morphism of complexes
$$
\lambda=i^*\oplus(i\sc)^*\fcolon C_G\updot\to C_T\updot\;,
$$
where $C_G\updot$ and $C_T\updot$ are the complexes
introduced in~\ref{subsec:morphisms-of-complexes}:
\begin{gather*}
C_G\updot=[\KD(\Gbar)_1\to \KD(\Gbar\sc)_1\rangle,  \\
C_T\updot=[\KD(T)_1\to \KD(T\sc)_1\rangle.
\end{gather*}
In other words,
\begin{gather*}
C_G\updot=\KD(\Gbar)_1\oplus \KD(\Gbar\sc)_1[-1],  \\
C_T\updot=\KD(T)_1\oplus \KD(T\sc)_1[-1],
\end{gather*}
with  differentials given by the matrix of formula \eqref{eq:matrix-fibre}
in \ref{subsec:cones-and-fibres}.

To prove the theorem, it suffices to prove that $\lambda$ is a quasi-isomorphism.
We must prove that the maps
$$
\lambda^0\fcolon \sH^0(C_G\updot)\to \sH^0(C_T\updot)
\text{ and }
\lambda^1\fcolon \sH^1(C_G\updot)\to \sH^1(C_T\updot)
$$
are isomorphisms.

We prove that $\lambda^0$ is an isomorphism.
Using Rosenlicht's lemma, we see immediately that
\begin{gather*}
\sH^0(C_G\updot)=\ker[\rho^*\fcolon \Chi^*(\Gbar)\to\Chi^*(\Gbar\sc)],\\
\sH^0(C_T\updot)=\ker[\rho_T^*\fcolon \Chi^*(\Tbar)\to\Chi^*(\Tbar\sc)],
\end{gather*}
and $\lambda^0\fcolon \sH^0(C_G\updot)\to \sH^0(C_T\updot)$ is the map
induced by the restriction map $i^*\fcolon \Chi^*(\Gbar)\to \Chi^*(\Tbar)$.
Now it is clear that $\lambda^0$ is an isomorphism.

We prove that $\lambda^1$ is an isomorphism.
Write $Z=\ker\:\rho$.
Consider the composed map
$$
\sigma_G\fcolon \Chi^*(\Zbar)\labelto{c} \Pic(\Gbar\sss)\isoto
\sH^1(\KDiv(\Gbar)_1)\isoto\sH^1(C\updot_G).
$$
Here the last isomorphism is induced by the isomorphism
$\KDiv(\Gbar)_1\isoto C\updot_G$ in the derived category
and the map $c$ is the map of Construction \ref{constr:FI}.
The map $\sigma_G$ is an isomorphism, because it is a composition
of isomorphisms.
We compute $\sigma_G$ explicitly.

We have
\begin{gather*}
C_G^1=\Div(\Gbar)_e\oplus \kbar(\Gbar\sc)^\times_{e,1}  \\
\ker d_G^1=\{(D_G, f_\Gsc)\in \Div(\Gbar)_e\oplus \kbar(\Gbar\sc)^\times_{e,1}
\ |\ \rho^*(D_G)=\divisor( f_\Gsc)\}
\end{gather*}
where $d^1_G\fcolon C^1_G\to C^2_G$ is the  differential in $C_G\updot$.

Let $\chi\in\Chi^*(\Zbar)$.
Define a right action of $\Zbar$ on $\Gbar\sc\times\Gabar$ by $(g,a)*z=( gz,\chi(z)^{-1}a)$,
where $g\in \Gbar\sc, a\in \Gabar, z\in \Zbar$.
Set $E\sss=L(\chi)=(\Gbar\sc\times \Gabar)/\Zbar$,
then $E\sss$ is a linear bundle over $\Gbar\sss=\Gbar\sc/\Zbar$.
By definition $\cl(E\sss)=c(\chi)\in\Pic(\Gbar\sss)$.
Since we have a canonical isomorphism $\Pic(\Gbar)\isoto \Pic(\Gbar\sss)$, our $E\sss$
comes from a unique (up to an isomorphism) line bundle $E$ over $\Gbar$.
Let $\varphi$ be a rational section of $E$ such that $\varphi(e)\neq 0,\infty$.
Set $D_G=\divisor(\varphi)$, then $\cl(D_G)$ is the image of $\chi$
in $\Pic(\Gbar)=\sH^1(\KDiv(\Gbar)_1)$.
Set $D_\Gsc=\rho^*(D_G)$.
Since $\Pic(\Gbar\sc)=0$, there exists $ f_\Gsc$ on $\Gbar\sc$ such that
$D_\Gsc=\divisor( f_\Gsc)$.
Since $e_G\notin\textup{supp}(D_G)$, we see that $f_\Gsc(e_\Gsc)\neq 0,\infty$.
Set $f'_\Gsc=f_\Gsc/f_\Gsc(e_\Gsc)$.
Then $(D_G, f'_\Gsc)\in\ker\: d_G^1$ and
$\cl(D_G,  f'_\Gsc)=\sigma_G(\chi)\in\sH^1(C_G\updot)$.

We need a lemma.
\end{subsec}

\begin{lemma}\label{lem:Z}
 The restriction of $f'_{\Gsc}$ to $Z$ is $\chi^{-1}$.
\end{lemma}

\begin{proof}[Proof of Lemma \ref{lem:Z}]

Consider the section $\rho^*(\varphi)$ of $\rho^* E$.
By the construction of $E\sss$ we have  a canonical trivialization
$$
\mu\fcolon  \Gbar\sc\times\Gabar \isoto  \rho^*E
$$
which maps $\rho^*(\varphi)$ to some $\psi=\mu^*(\rho^*(\varphi))$.
We have
\begin{equation}\label{eq:G-psi}
\psi(gz)=\chi(z)^{-1}\psi(g) \text{ for all }  g\in \Gbar\sc,\ z\in \Zbar
\end{equation}
because $\varphi|_{\Gbar\sss}$ is a rational section of $E\sss$.
But
$$
D_\Gsc=\rho^*(D_G)=\divisor(\rho^*(\varphi))=\divisor(\psi),
$$
so we may take $f_\Gsc=\psi$.
Since $U(\Gbar\sc)=\Chi^*(\Gbar\sc)=0$, there exists, up to a constant,
only one rational function $f_\Gsc$ on $\Gbar\sc$ such that $D_\Gsc=\divisor( f_\Gsc)$.
Using \eqref{eq:G-psi}, we obtain that for any such $f_\Gsc$ we have
$$
f'_\Gsc(z)=f_\Gsc(z)/f_\Gsc(e)=\psi(z)/\psi(e)=\chi(z)^{-1}.
$$
\end{proof}

\begin{subsec}
{\it Proof of Theorem \ref{thm:main} (cont.)}.
We have
\begin{gather*}
C_T^1=\Div(\Tbar)_e\oplus \kbar(\Tbar\sc)_{e,1}^\times,  \\
\ker d_T^1=\{(D_T, f_\Tsc)|\ \rho^*_T(D_T)=\divisor( f_\Tsc)\},
\end{gather*}
where $d^1_T\fcolon C^1_T\to C^2_T$ is the differential in $C\updot_T$.
The canonical isomorphism
\begin{equation}\label{eq:T-complex-iso}
C_T\updot\isoto [\Chi^*(\Tbar)\to\Chi^*(\Tbar\sc) \rangle
\end{equation}
induces a composed map
$$
\tau_T\fcolon \sH^1(C_T\updot)\isoto \sH^1([\Chi^*(\Tbar)
\to\Chi^*(\Tbar\sc) \rangle)
\isoto \Chi^*(\Zbar)
$$
where the latter isomorphism is defined
as follows: $\cl(\kappa)\mapsto\kappa|_\Zbar$
for  $\kappa\in\Chi^*(\Tbar\sc)$.
We compute $\tau_T$ explicitly.

Let $(D_T, f_\Tsc)\in\ker\: d_T^1$.
Since $\Pic(\Tbar)=0$, there exists a rational function
$f_T\in \kbar(\Tbar)_{e,1}^\times$
such that $\divisor(f_T)=-D_T$.
Set $\tilde{f}_\Tsc=f_\Tsc \cdot\rho^*(f_T)\in \kbar(\Tbar\sc)_{e,1}^\times$.
Then $\divisor(\tilde{f}_\Tsc)=0$ and $\tilde{f}_\Tsc(e)=1$.
By Rosenlicht's lemma $\tilde{f}_\Tsc\in\Chi^*(\Tbar\sc)$.
Moreover
$(0, \tilde{f}_\Tsc) \in \ker d_T^1$ and
$(0, \tilde{f}_\Tsc) = (D_T,f_\Tsc) + d^0_T(f_T)$.
The construction of the isomorphism~\eqref{eq:T-complex-iso}
then implies that
$\cl(D_T,f_\Tsc) \in \sH^1(C_T\updot)$ corresponds to
$\cl(\tilde{f}_\Tsc)\in \sH^1([\Chi^*(\Tbar)\to\Chi^*(\Tbar\sc\rangle)$.
The image of $\cl(\tilde{f}_\Tsc)$ in $\Chi^*(\Zbar)$ is
$$
\tilde{f}_\Tsc |_Z= f_\Tsc\cdot \rho^*(f_T)|_Z= f_\Tsc |_Z\;.
$$
Thus the map $\tau_T$ is given by $\cl(D_T,f_\Tsc)\mapsto f_\Tsc |_Z$.

Now we see that the composed map
$$
\beta\fcolon \Chi^*(\Zbar)\labelto{\sigma_G}\sH^1(C_G\updot)\labelto{\lambda^1}\sH^1(C_T\updot)
\labelto{\tau_T} \Chi^*(Z)
$$
is given by
$$
\chi\mapsto (D_G, f'_\Gsc)\mapsto (i^*(D_G), (i\sc)^*( f'_\Gsc))
\mapsto (i\sc)^*(f'_\Gsc)|_Z\;.
$$
Clearly we have $(i\sc)^*(f'_\Gsc)|_Z=f'_\Gsc|_Z$.
By Lemma \ref{lem:Z} $f'_\Gsc|_Z=-\chi$ (with additive notation).
We see that our composed  map $\beta$ is given by $\chi\mapsto -\chi$,
hence it is an isomorphism.
Since $\sigma_G$ and $\tau_T$ are isomorphisms, we conclude that $\lambda^1$
is an isomorphism.  This completes the proof of  Theorem \ref{thm:main}.\qed
\end{subsec}
\begin{remark}
A different proof of the existence of an isomorphism $\UPic(\Gbar)\isoto\pi_1(\Gbar)^D$
was proposed in \cite{CT06}.
\end{remark}

\section{UPic of torsors and the elementary obstruction}
\label{sec:upic-torsors-elt-obstr}

\begin{lemma}\label{lem:additivity}
Let $X$ and $Y$ be smooth geometrically integral $k$-varieties.
Assume that $Y$ is $\kbar$-rational.
Then the canonical morphism
$$
\zeta\colon\UPic(\Xbar)\oplus \UPic(\Ybar)\to\UPic(\Xbar\times\Ybar)
$$
induced by the projections $p_X$ and $p_Y$ from $X \times Y$
to the corresponding factors, is  a quasi-isomorphism.
\end{lemma}

\begin{proof}
By Rosenlicht's lemma \cite[Lemma 2.1]{Fossum-Iversen} the map
$$
\sH^0(\zeta)\colon U(\Xbar)\oplus U(\Ybar)\to U(\Xbar\times\Ybar)
$$
is an isomorphism. By a lemma of Colliot-Th\'el\`ene and Sansuc
\cite[Lemme 11 p. 188]{CTS77} the map
$$
\sH^1(\zeta)\colon\Pic(\Xbar)\oplus \Pic(\Ybar)\to\Pic(\Xbar\times\Ybar)
$$
is an isomorphism.
Thus $\zeta$ is a quasi-isomorphism.
\end{proof}

For a  $k$-torsor $X$ under a connected linear algebraic $k$-group $G$
it was shown by Sansuc that $U(\Xbar) = U(\Gbar) = \Chi^*(\Gbar)$ and
$\Pic(\Xbar) = \Pic(\Gbar)$.
Sansuc's result extends to $\UPic$, and so does his proof.

\begin{lemma}\label{lem:upic-torsor}
Let $\varsigma\fcolon X\times G\to X$ be a $k$-morphism defining a right action
of a connected linear algebraic $k$-group $G$ on a  smooth geometrically integral $k$-variety $X$.
Then
\begin{theoremlist}
\item
The canonical morphism
$$
\zeta\fcolon\UPic(\Xbar) \oplus \UPic(\Gbar)\to
\UPic(\Xbar \times \Gbar)
$$
is  a quasi-isomorphism.
\item
Denote by
\begin{equation*}
 \pi_G\fcolon \UPic(\Xbar \times \Gbar) =
\UPic(\Xbar) \oplus \UPic(\Gbar) \to \UPic(\Gbar)
\end{equation*}
the projection.
Then
$$
\varphi=\pi_G\circ\varsigma^*\fcolon \UPic(\Xbar)\to\UPic(\Xbar \times \Gbar)\to\UPic(\Gbar)
$$
is a canonical morphism, functorial in $(X,G)$ and equal, for any
$x_0\in X(k)$ to $\alpha^*_{x_0}$, where $\alpha_{x_0}\fcolon G\to X$ is the
$k$-morphism defined by $\alpha_{x_0}(g)=x_0 g$ for $g\in G$.
\item
If in addition $X$ is a torsor of $G$ over $k$, then $\varphi$
is an isomorphism in the derived category.
\end{theoremlist}
\end{lemma}

\begin{proof}
(i) Since $G$ is $\kbar$-rational, by Lemma \ref{lem:additivity} $\zeta$ is a quasi-isomorphism.

(ii) Take $x_0 \in X(k)$. Let $i_{x_0}$ be the $k$-morphism
$G \to X\times G$ defined by $i_{x_0}(g) = (x_0, g)$.
Then $p_G \circ i_{x_0} = \id$ and $p_X \circ i_{x_0}$ is the
constant map $G \to \{x_0\} \in X$.
Hence $i_{x_0}^* = \pi_G$ and, since $\alpha_{x_0} = \varsigma \circ i_{x_0}$, we get
\[\alpha_{x_0}^* = i_{x_0}^* \circ \varsigma^* = \pi_G
\circ \varsigma^* = \varphi. \]

(iii) By  \cite[Lemmes 6.4, 6.5(ii), 6.6(i)]{Sansuc:brauer-gal}
the morphisms $\sH^0(\varphi)\fcolon U(\Xbar)\to\Chi^*(\Gbar)$
and $\sH^1(\varphi)\fcolon \Pic(\Xbar)\to\Pic(\Gbar)$
are isomorphisms,
hence $\varphi$ is an isomorphism in the derived category.
\end{proof}

As a corollary we obtain a canonical isomorphism between
the target of the elementary
obstruction (Definition \ref{def:elementary-obstruction})
and  the abelian Galois
cohomology  $H^i_\ab(k,G):=H^i(k,\langle T\sc\to T])$.

\begin{corollary}\label{cor:elt-obstruction-target}
Let $X$ be a $k$-torsor under a connected linear algebraic $k$-group $G$.
We have a canonical isomorphism
\[
 \Ext^i(\UPic(\Xbar), \Gm) \simeq H^i_\ab(k, G)
\]
which is functorial in $(G, X)$ and in~$k$.
\end{corollary}

\begin{proof}
By Lemma \ref{lem:upic-torsor}(iii) $\Ext^i(\UPic(\Xbar),\Gm)=\Ext^i(\UPic(\Gbar),\Gm)$.
By Corollary \ref{cor:upic-dual-tori} $\Ext^i(\UPic(\Gbar),\Gm)=H^i_\ab(k,G)$.
\end{proof}

Let $\ab^1\fcolon H^1(k,G)\to H^1_\ab(k,G)$ be the abelianization map
constructed in \cite{Borovoi:Memoir}.
For a $k$-torsor $X$ of $G$, let $[X]$ denote its class in $H^1(k,G)$.
We write $[X]_\ab:=\ab^1([X])\in H^1_\ab(k,G)$.
We shall prove that
the elementary obstruction $e(X)$ coincides up to sign with the  $[X]_\ab$.
For semisimple groups this was proved by Skorobogatov
\cite[p. 54]{Skorobogatov:torsors}.
For tori this was proved by Sansuc \cite[(6.7.3), (6.7.4)]{Sansuc:brauer-gal}
(see also Skorobogatov \cite[Lemma 2.4.3]{Skorobogatov:torsors}).
First we give Sansuc's proof for tori with details added.

\begin{lemma}[Sansuc]\label{lem:torus-elt-obstr}
Let $T$ be a torus over $k$ and let $X$ be a $k$-torsor under $T$,
determined by a cocycle  $c\fcolon \sigma\mapsto c_\sigma\fcolon \Gal(\kbar/k)\to T(\kbar)$.
Consider the extension
\begin{equation}
1\to\kbar^\times\to\kbar[X]^\times\to \Chi^*(\Tbar)\to 1
\end{equation}
The class $e(X)$ of this extension in
$\Ext^1(\Chi^*(\Tbar), \kbar^\times)$ corresponds
 under the isomorphism of Lemma \ref{lem:Ext-H^i}
\[
\Ext^1(\Chi^*(\Tbar), \kbar^\times) = H^1(k,T)
\]
to the class of the cocycle $c^{-1}$.
\end{lemma}

\begin{proof}
We regard $X(\kbar)$ as $T(\kbar)$ with the twisted Galois action
$\sigma *t=c_\sigma \cdot \sigma t$, where $t\in T(\kbar)$. Similarly we regard
$\kbar[X]^\times$ as $\kbar[T]^\times$ with the twisted Galois action, etc.
In all cases we use the notation
$\sigma*$ to denote the twisted Galois action.

Let $\chi\in\Chi^*(\Tbar)$.
We compute ${}^{\sigma*}\chi$.
For $t\in T(\kbar)$ we have
$$
({}^{\sigma*}\chi)(\sigma*t)=\sigma(\chi(t)),
$$
hence
\begin{multline*}
({}^{\sigma*}\chi)(t)= \sigma(\chi(\sigma^{-1}*t))
=\sigma(\chi(c_{\sigma^{-1}}\cdot \sigma^{-1}t)) =
\\
=\sigma(\chi(\sigma^{-1}(c_\sigma^{-1} t)))
= {}^\sigma\chi(c_\sigma^{-1} t)
={}^\sigma\chi(c_\sigma^{-1})\cdot {}^\sigma\chi(t).
\end{multline*}
Thus
$$
{}^{\sigma*}\chi={}^{\sigma}\chi(c_\sigma^{-1})\cdot {}^{\sigma}\chi.
$$

Now let $\varphi\fcolon \Chi^*(\Tbar)\to\kbar[\Xbar]^\times$ be the
standard (non-equivariant) splitting corresponding to the
identification of $\Xbar$ with $\Tbar$.
By abuse of notation we denote this splitting by $\chi\mapsto \chi$.
Then the extension class
$e(X) \in H^1(k, \Hom_\Zz(\Chi^*(\Tbar), \kbar^\times))$
is represented by the cocycle
$\sigma\mapsto \sigma \varphi\cdot \varphi^{-1}$.

Since
$$
(\sigma\varphi)(\chi)={}^{\sigma*}(\varphi({}^{\sigma^{-1}}\chi))
={}^{\sigma*}({}^{\sigma^{-1}}\chi)
=\chi(c_\sigma^{-1})\cdot \chi,
$$
we see that $e(X)$ is represented by the cocycle
\[
\sigma \mapsto \sigma\varphi\cdot\varphi^{-1} =
(\chi\mapsto \chi(c_\sigma^{-1})) \in \Hom_\Zz(\Chi^*(\Tbar), \kbar^\times),
\]
which corresponds to the cocycle $\sigma \mapsto c_\sigma^{-1}$ under
the identification $T(\kbar) \simeq \Hom_\Zz(\Chi^*(\Tbar), \kbar^\times)$
given by $t\mapsto\, (\chi\mapsto\chi(t))$ for $t\in T(\kbar)$.
\end{proof}

\begin{theorem}\label{th:elementary-abcoh}
Let $X$ be a torsor under a linear algebraic group $G$ over $k$. The
elementary obstruction class $e(X) \in \Ext^1(\UPic(\Xbar), \kbar^\times)$
corresponds to  $-[X]_\ab \in H^1_\ab(k, G)$
under the canonical isomorphism of Corollary~\ref{cor:elt-obstruction-target}
\[
 \Ext^i(\UPic(\Xbar), \Gm) \simeq H^i_\ab(k, G).
\]
\end{theorem}
\begin{proof}
Without loss of generality we may assume that $G$ is reductive.

As in \cite[p.~369]{Kot86}
(compare \cite[Lemma 1.1.4(i)]{Borovoi-Kunyavskii:formulas}),
we construct an epimorphism $\alpha\fcolon H\to G$, where $H$
is a reductive $k$-group with $H\sss$ simply connected,
together with a $k$-torsor $X_H$ under $H$
such that $\alpha_*(X_H)\simeq X$.
By functoriality of the isomorphism of Corollary~\ref{cor:elt-obstruction-target},
in order to prove the theorem for $G$ and $X$,   it is sufficient to prove it
for $H$ and $X_H$.

Since $H\sss$ is simply connected, the homomorphism $H \to H\tor$ induces an isomorphism
$H^1_\ab(k, H) \simeq H^1_\ab(k, H\tor)= H^1(k,H\tor)$, cf. \cite[Example 2.12(2)]{Borovoi:Memoir}.
We see that the functoriality of the isomorphism
of Corollary~\ref{cor:elt-obstruction-target}
implies that it is sufficient to prove the theorem for torsors under
tori, which was done in Lemma~\ref{lem:torus-elt-obstr}.
\end{proof}

\begin{corollary}\label{cor:elt-obstruction-cpct}
  Let $Y$ be a smooth compactification of a torsor $X$ under a connected linear
algebraic group $G$ over $k$.
Let $S$ be the N\'eron-Severi torus of $Y$, i.e. the $k$-torus $S$ such that
$\Chi^*(S)=\Pic(\Ybar)$.
Then $\Ext^1(\UPic(\Ybar),\Gm)=H^2(k,S)$ and
we have a canonical injection
\begin{equation*}
H^1_\ab(k, G)\into H^2(k, S)
\end{equation*}
sending  $[X]_\ab\in H^1_\ab(k, G)$ to
 $-e(Y) \in H^2 (k, S)$.
\end{corollary}
\begin{proof}
By Corollary \ref{cor:restr-elementary} the open embedding $j\fcolon X\into Y$ induces
an injection $j_* \fcolon \Ext^1(\UPic(\Xbar),\Gm)\into \Ext^1(\UPic(\Ybar),\Gm)$.
By Corollary \ref{cor:elt-obstruction-target} $\Ext^1(\UPic(\Xbar),\Gm)=H^1_\ab(k,G)$.
Since $\UPic(\Ybar)=\Pic(\Ybar)[-1]$, we have
$\Ext^1(\UPic(\Ybar),\Gm)=\Ext^2(\Pic(\Ybar),\Gm)=H^2(k,S)$
(we use Lemma \ref{lem:Ext-H^i}).
By Theorem \ref{th:elementary-abcoh} $j_*$ takes $[X]_\ab$ to $-e(Y)$.
\end{proof}

\begin{proposition}  \label{prop:torsor-p-adic}
For (a smooth compactification of) a torsor $X$
under a connected linear algebraic
group $G$ over a $p$-adic field $k$, the
elementary obstruction is the only obstruction to the existence of
$k$-rational points.
\end{proposition}

\begin{proof}
We will first show that the vanishing of $e(X)$ implies the existence
of a $k$-rational point on $X$.
Let $\ab^1\fcolon H^1(k,G)\to H^1_\ab(k,G)$
denote the abelianization map of \cite{Borovoi:Memoir}.
We have an exact sequence
\[
H^1(k, G\sc)\to H^1(k, G) \labelto{\ab^1} H^1_\ab(k, G),
\]
see~\cite[(3.10.1)]{Borovoi:Memoir}.
By Theorem~\ref{th:elementary-abcoh} $\ab^1([X])=-e(X)$,
and by assumption $e(X)=0$.
By Kneser's theorem $H^1(k,G\sc)=0$.
We conclude that $[X]=0$.
Thus $X$ has a $k$-point.

Now we will show that for
a smooth compactification $Y$ of $X$
the vanishing of
$e(Y)$ implies the existence of a $k$-rational point.
By Corollary~\ref{cor:restr-elementary} the vanishing of $e(Y)$
implies the vanishing of $e(X)$.
As above, this implies that $X(k)\neq\emptyset$, hence $Y(k)\neq\emptyset$.
\end{proof}

\begin{proposition}\label{prop:torsor-number}
For (a smooth compactification of) a torsor $X$ under a connected linear algebraic
group $G$ over a number field $k$, the
elementary obstruction is the only obstruction to the Hasse
principle.
\end{proposition}
\begin{proof}
First assume that $X(k_v)\neq\emptyset$ for all places $v$ of $k$,
and assume that $e(X)=0$.
Clearly $[X]\in\Sha^1(k,G)$, where $\Sha^1(k,G)$ is the Tate-Shafarevich kernel for $G$.
By Theorem~\ref{th:elementary-abcoh} $\ab^1([X])=-e(X)$.
Clearly $e(X)\in\Sha^1_\ab(k,G)$, where
$$
\Sha^1_\ab(k,G):=\ker\left[H^1_\ab(k,G)\to\prod_v H^1_\ab(k_v, G)\right].
$$
By \cite[Thm.~5.12]{Borovoi:Memoir} the induced map
$\ab^1_\Sha\fcolon\Sha^1(k, G) \to \Sha^1_\ab(k, G)$ is bijective
(here the Hasse principle for semisimple simply connected
groups, due to Kneser, Harder and Chernousov, plays  a major
role in the proof).
We see that $[X]=0$, hence $X(k)\neq\emptyset$.

Now let $Y$ be a smooth compactification of a $k$-torsor $X$.
Assume that $e(Y)=0$ and that $Y(k_v)\neq\emptyset$ for all $v$.
By Corollary~\ref{cor:restr-elementary} $e(X)=0$ because $e(Y)=0$.
Since $Y$ is smooth, we have $X(k_v)\neq\emptyset$ for all $v$.
As above we see that $X(k)\neq\emptyset$, hence $Y(k)\neq\emptyset$.
\end{proof}

For  other proofs of Propositions \ref{prop:torsor-p-adic} and \ref{prop:torsor-number}
see \cite{BCTS06}.

\providecommand{\url}[1]{\textit{#1}}
\providecommand{\urlbr}[1]{\discretionary{}{}{}}

\providecommand{\bysame}{\leavevmode\hbox to3em{\hrulefill}\thinspace}

\end{document}